\documentclass[10pt]{amsart}
\usepackage{amsmath}
\usepackage{amssymb,amsgen,amsbsy,amsopn,amsfonts,graphicx,amsthm}
\usepackage[dvips]{color}
\usepackage{multicol}
\usepackage{mathrsfs}

     \font\tenmsb=msbm10
     \font\sevenmsb=msbm7
     \font\fivemsb=msbm5
     \catcode`\@=11
     \ifx\amstexloaded@\relax\catcode`\@=\active
     \endinput\else\let\amstexloaded@\relax\fi
     \def\spaces@{\space\space\space\space\space}
     \def\spaces@@{\spaces@\spaces@\spaces@\spaces@\spaces@}
     \def\space@.{\futurelet\space@\relax}
     \space@. %
     \def\Err@#1{\errhelp\defaulthelp@\errmessage{AmS-TeX error: #1}}
     \def\relaxnext@{\let\next\relax}
     \def\accentfam@{7}
     \def\noaccents@{\def\accentfam@{0}}
     \def\Cal{\relaxnext@\ifmmode\let\next\Cal@\else
     \def\next{\Err@{Use \string\Cal\space only in math mode}}\fi\next}
     \def\Cal@#1{{\Cal@@{#1}}}
     \def\Cal@@#1{\noaccents@\fam\tw@#1}
     \def\Bbb{\relaxnext@\ifmmode\let\next\Bbb@\else
     \def\next{\Err@{Use \string\Bbb\space only in math mode}}\fi\next}
     \def\Bbb@#1{{\Bbb@@{#1}}}
     \def\Bbb@@#1{\noaccents@\fam\msbfam#1}
     \newfam\msbfam
     \textfont\msbfam=\tenmsb
     \scriptfont\msbfam=\sevenmsb
     \scriptscriptfont\msbfam=\fivemsb

\newtheorem{theorem}{Theorem}[section]

\newtheorem{definition}{Definition}[section]

\newtheorem{lemma}{Lemma}[section]

\numberwithin{equation}{section}
\newtheorem{remark}{Remark}[section]

\begin{document}

\title[Analysis of solid tumor growth]{Analysis of a mixture model of tumor growth}
\author{John Lowengrub}
\address{Department of Mathematics, and Department of Biomedical Engineering, University of California, Irvine, CA 92697-3875, USA.}
\email{lowengrb@math.uci.edu.}
\author{Edriss S. Titi}
\address{Department of Mathematics, and Department of Mechanical and Aerospace Engineering, University of California, Irvine, CA 92697-3875, USA. Also, The Department of
Computer Science and Applied Mathematics, The Weizmann Institute of Science, Rehovot 76100, Israel. Fellow of the Center of Smart Interfaces (CSI), Technische Universit\"{a}t Darmstadt, Germany.}
\email{etiti@math.uci.edu.}
\author{Kun Zhao}
\address{Department of Mathematics, University of Iowa, Iowa City, IA 52242-1419, USA.}
\email{kun-zhao@uiowa.edu.}


\begin{abstract}
We study an initial-boundary value problem (IBVP) for a coupled Cahn-Hilliard-Hele-Shaw system that models tumor growth. For large initial data with finite energy, we prove global (local resp.) existence, uniqueness, higher order spatial regularity and Gevrey spatial regularity of strong solutions to the IBVP in 2D (3D resp.). Asymptotically in time, we show that the solution converges to a constant state exponentially fast as time tends to infinity under certain assumptions.
\end{abstract}

\maketitle

\section{Introduction}
Approximately 550,000 patients will die from cancer before the ball drops in New York City's Time's Square again. In the United Kingdom, one in four people will die of cancer, whilst one in three will at some point in their life be diagnosed to have cancer (http://www.cancerresearchuk.com/). In comparison to molecular biology, cell biology, and drug delivery research, mathematics has so far contributed relatively little to the area. A search in the PubMed bibliographic database (http://www.ncbi.nlm.nih.gov/PubMed/) shows that out of 1.5 million papers in the area of cancer research, approximately 5$\%$ are concerned with mathematical modeling. However, it is clear that mathematics could make a huge contribution to many areas of experimental cancer investigation since there is now a wealth of experimental data which requires systematic analysis. At the current stage of cancer research, most of the mathematical models are built and developed from three perspectives: discrete (microscopic), continuous (macroscopic) and hybrid (micro-macroscopic). This paper focuses on a continuous mixture model. The model is designed to capture the dynamics of morphological changes in solid tumor growth. The thermodynamically consistent model is derived using a classical continuum mechanics approach, based on the principle of mass conservation and the second law of thermodynamics. The free energy is chosen by taking into account the effect of cell-cell and cell-matrix adhesion.
In this approach, sharp tumor/host interfaces are replaced by narrow transition layers. These models are capable of describing the avascular, vascular and metastasis stages of  solid tumor growth. In non-dimensional form, the simplest version of the equations reads (c.f. \cite{FJCWLC,WLFC}):
\begin{equation}
\left\{
\begin{aligned}
&\phi_t+\vec{\mathbf{v}}\cdot\nabla\phi=\Delta \mu+S,\ \ \mathbf{x}\in\Omega,\ \ t>0,\\
&\mu:=F'(\phi)-\epsilon^2\Delta\phi=\phi^3-\phi-\epsilon^2\Delta\phi,\\
&\vec{\mathbf{v}}=-\nabla P+\gamma\mu\nabla\phi,\\
&\nabla\cdot \vec{\mathbf{v}}=S.
\end{aligned}
\right.\label{e11}
\end{equation}
Here, $\phi$ is the scalar volume fraction (order parameter) of tumor cells, $\mu$ is the chemical potential, $S$ is a mass source term that accounts for cell proliferation, $\vec{\mathbf{v}}$ denotes the advective vector velocity field, $P$ is the scalar pressure, and $\epsilon>0,\gamma>0$ are constants with the former describing the interface thickness and the latter describing cell-cell adhesion. The function $F(\phi)$, which appears in the definition of the chemical potential, i.e. in the second equation of (\ref{e11}) takes the form $F(\phi)=\frac14(\phi^2-1)^2$, which has a physically relevant, double-well structure, each of them representing the two phases of the mixture. We observe that, when $S=0$, system (\ref{e11}) reduces to a simplified version of an earlier model derived in \cite{GLL,GLL-1}, that is used to describe gravity-driven, density-mismatched, two-phase flow in a Hele-Shaw cell. A similar set of equations was also used in \cite{OS} as a model for spinodal decomposition of a binary fluid in a Hele-Shaw cell. In this paper, we are interested in the case in which $\gamma$ and $\epsilon$ are strictly positive. Therefore, without loss of generality, we assume that $\epsilon=\gamma=1$ throughout the paper.

We consider system (\ref{e11}) in a bounded domain $\Omega$ which is a rectangle in $\mathbb{R}^2$ or a box in $\mathbb{R}^3$. The system is supplemented by the following initial-boundary value conditions:
\begin{equation}
\left\{
\begin{aligned}
&(\phi,\mu)(\mathbf{x},0)=(\phi^0,\mu^0)(\mathbf{x}),\ \ \mathbf{x}\in\Omega;\\
&\vec{\mathbf{v}}\cdot\mathbf{n}|_{\partial\Omega}=0,\\
&\nabla \phi\cdot\mathbf{n}|_{\partial\Omega}=\nabla \mu\cdot\mathbf{n}|_{\partial\Omega}=0,\ t\ge0,
\end{aligned}
\right.\label{e12}
\end{equation}
where $\Omega=[0,L]\times[0,l]$ in $\mathbb{R}^2$ and $\Omega=[0,L]\times[0,l]\times[0,h]$ in $\mathbb{R}^3$, and $\mathbf{n}$ is the
unit outward normal to $\partial\Omega$. In view of the definition of $\mu$ and the boundary conditions (\ref{e12}) we see that, in fact, the boundary condition for $\mu$ is equivalent to $\nabla(\Delta\phi)\cdot\mathbf{n}|_{\partial\Omega}=0$.

System (\ref{e11}) is closely related to the Cahn-Hilliard-Navier-Stokes type system of equations studied in \cite{ABP,BF,Boyer,Boyer-1,Boyer-2,DD,EG,GPV,KW,KKL,M,WW} and many others. In the absence of the source term, system (\ref{e11}) serves as a model for cell sorting owing to differential cell-to-cell adhesion \cite{WLFC}.

The CHHS system (\ref{e11}), with $S=0$, was studied numerically in \cite{W}, where a convex-splitting scheme was proposed and was proven to be unconditionally stable and solvable. A practical and efficient multigrid method was used for solving the scheme at each time step. However, to the best of our knowledge, system (\ref{e11}) has not been studied analytically, even when $S=0$. In this paper, we will study qualitative behavior of solutions to (\ref{e11})--(\ref{e12}) with $S=0$, i.e., existence, uniqueness, regularity and long-time behavior of solutions to (\ref{e11})--(\ref{e12}) with $S=0$. Our first goal is to prove well-posedness of strong solutions, see Definition 2.1 below, to the IBVP. We show that for initial data belonging to $H^2$, there exists a unique global (local resp.) strong solution to (\ref{e11})--(\ref{e12}) of the 2D (3D resp.) system. The second goal of this paper is to improve the spatial regularity of the strong solution for the same initial data. Following ideas used for parabolic PDEs and in particular for the Navier-Stokes equations (c.f. \cite{CF,T}), we show that the spatial regularity of the solution increases by degree 2 within the lifespan of the solution. Thirdly, following \cite{FT,LT} (see also \cite{FT-2}), we study the Gevrey regularity (spatial analyticity) of the solution in the domain $\Omega$. We show that for initial data belonging to $H^4$, the solution belongs to the Gevrey space $G_t^2(\Omega)$, see definition in Section 8, globally in 2D and locally in 3D. In other words, the solution is globally (locally resp.) spatially analytic in 2D (3D resp.). Since the elements in the Gevrey class have high-mode coefficients that decay exponentially in wave number to zero, as a result of the Gevrey regularity one can show that the standard Galerkin method, based on basis functions in (\ref{e21})--(\ref{e23}), converges exponentially fast, see e.g. \cite{DT1,DT2,GST,JMT,LT}. The last part of this paper is devoted to the study of the long-time dynamics of solutions to (\ref{e11})--(\ref{e12}). First, we show that $\phi$ converges exponentially rapidly to its spatial average $\bar{\phi}$ over the domain $\Omega$ in the $H^1$ norm, as time goes to infinity, provided that $\bar{\phi}$ lies outside the spinodal region where $F''(\cdot)=3\phi^2-1\ge0$ and the initial perturbation $\|\phi^0-\bar{\phi}\|^2_{H^2}$ is sufficiently small. The result holds true in both 2D and 3D. It should be pointed out that the idea for our proof is in the spirit of \cite{Boyer}, where the author studied the long-time behavior of strong solutions to an initial-boundary value problem for a coupled Cahn-Hilliard-Navier-Stokes system. On the other hand, we show that, in the 2D case, without the smallness assumption on the initial perturbation, $\phi$ still converges exponentially fast to $\bar{\phi}$ in the $H^2$ norm, as time goes to infinity, provided that $l\le L<\pi$. This condition indeed implies that the constant $C_\Omega$ in Poincar\'{e}'s inequality (i.e., $\|f-\bar{f}\|_{L^2}^2\le C_\Omega \|\nabla f\|^2_{L^2}$) on the domain $\Omega=[0,L]\times[0,l]$ satisfies $C_\Omega<1$ since, for the domain under consideration, it can be explicitly computed that $C_\Omega=(\lambda_1)^{-1}=L^2/\pi^2$ using the eigenfunctions (\ref{e21}) given below, where $\lambda_1$ is the first eigenvalue of the Laplacian operator subject to the Neumann boundary condition. Regarding the condition $l\le L<\pi$ (or $C_\Omega<1$), we remark that, since the interface thickness $\epsilon$ is assumed to be $1$ in the beginning, the units are still consistent. Indeed, in terms of $\epsilon$, the condition is $l\le L<\epsilon\pi$ for the original model \eqref{e11}.

Here, throughout this paper, $\|\cdot\|_{L^p}$, $\|\cdot\|_{L^\infty}$
and $\|\cdot\|_{W^{s,p}}$ denote the norm of the usual Lebesgue
measurable function spaces $L^p$ ($1\le p<\infty$), $L^\infty$ and
the usual Sobolev space $W^{s,p}$, respectively. For $p=2$, we
denote the norm $\|\cdot\|_{L^2}$ by $\|\cdot\|$ and
$\|\cdot\|_{W^{s,2}}$ by $\|\cdot\|_{H^s}$, respectively. The function
spaces under consideration are $$C([0,T];H^r(\Omega))\ \
\mathrm{and}\ \ L^2([0,T];H^s(\Omega))$$ equipped with norms
$$
\sup_{0\le t\le T}\|\Psi(\cdot,t)\|_{H^r}\ \ \mathrm{and}\ \ \big(\int_0^T\|\Psi(\cdot,\tau)\|_{H^s}^2d\tau\big)^{1/2},
$$
where $r,s$ are positive integers.

Let
$$
\begin{aligned}
\mathcal{V}&:=\{\theta\in C^\infty(\overline{\Omega})\ |\ \nabla\theta\cdot\mathbf{n}|_{\partial\Omega}=\nabla(\Delta\theta)\cdot\mathbf{n}|_{\partial\Omega}=0\},
\end{aligned}
$$
and
$$
\begin{aligned}
H&:=\mathrm{Closure\ of\ }\mathcal{V}\ \mathrm{in\ the}\ L^2\ -\ \mathrm{norm},\\
V&:=\mathrm{Closure\ of\ }\mathcal{V}\ \mathrm{in\ the}\ H^2\ -\ \mathrm{norm},\\
W&:=\mathrm{Closure\ of\ }\mathcal{V}\ \mathrm{in\ the}\ H^4\ -\ \mathrm{norm}.
\end{aligned}
$$

Throughout the paper, unless specified, $C$ will
denote a generic constant which is independent of the unknown functions
$\phi,\mu$ and $\vec{\mathbf{v}}$, but may depend on $\Omega$, initial data and the time $T$.
The value of $C$ may vary line by line according to
the context.

Our main results, for \eqref{e11} with $S=0$, are stated in the following theorems. First, for the existence of strong solutions, we have

\begin{theorem}[2D Global existence]
Let $\Omega=[0,L]\times[0,l]$ and let $\phi^0\in V$ be given. Then for every $T>0$ there exists a global strong solution $(\phi,\mu,\vec{\mathbf{v}},P)$, see Definition 2.1 below, to the initial-boundary value problem (\ref{e11})--(\ref{e12}) such that $\phi\in C([0,T];V)\cap L^2([0,T];W)$, $\mu\in C([0,T];H)\cap L^2([0,T];V)$, $\vec{\mathbf{v}}\in L^2([0,T];(H^1(\Omega))^2)$ and $P\in L^2([0,T];H^2(\Omega))$.
\end{theorem}

\begin{theorem}[3D Local existence]
Let $\Omega=[0,L]\times[0,l]\times[0,h]$ and let $\phi^0\in V$ be given. Then there exists a local strong solution $(\phi,\mu,\vec{\mathbf{v}},P)$ to the initial-boundary value problem (\ref{e11})--(\ref{e12}) such that $\phi\in C([0,T^*];V)\cap L^2([0,T^*];W)$, $\mu\in C([0,T^*];H)\cap L^2([0,T^*];V)$, $\vec{\mathbf{v}}\in L^2([0,T^*];(H^1(\Omega))^3)$ and $P\in L^2([0,T^*];H^2(\Omega))$ for some finite $T^*>0$, which depends on the initial data and the domain $\Omega$.
\end{theorem}

Concerning the uniqueness and well-posedness of the strong solutions obtained in the above theorems, we have

\begin{theorem}[2D $\&$ 3D Uniqueness and well-posedness]
Let the conditions of Theorem 1.1 and Theorem 1.2 hold. Then the solution obtained in Theorem 1.1 (Theorem 1.2 resp.) is unique and depends continuously on the initial data.
\end{theorem}

For higher order spatial regularity of the solutions obtained in Theorem 1.1 and Theorem 1.2, we have

\begin{theorem}[2D $\&$ 3D Higher order spatial regularity]
Let the conditions of Theorem 1.1 and Theorem 1.2 be satisfied and let $\phi$ be the strong solution to (\ref{e11})--(\ref{e12}). Then
$$
\begin{aligned}
&\|\phi(\cdot,t)\|^2_{H^4}
\le \frac{C(\mathcal{T})}{t},\ \ \forall\ t\in(0,\mathcal{T}),\\
&\int_0^t \tau^{\alpha+1}\|\phi(\cdot,\tau)\|_{H^6}^2d\tau\le C(\mathcal{T})\frac{t^\alpha}{\alpha},\ \ \forall\ t\in(0,\mathcal{T}),
\end{aligned}
$$
where $C(\mathcal{T})$ is an increasing function of $\mathcal{T}$, $\alpha$ is an arbitrary positive constant, and $\mathcal{T}>0$ is the lifespan of the solution, i.e., $\mathcal{T}>0$ is arbitrary for the 2D case and $\mathcal{T}=T^*$ for the 3D case.
\end{theorem}

The next theorem gives the Gevrey regularity of the solution to (\ref{e11})--(\ref{e12}).

\begin{theorem}[2D $\&$ 3D Gevrey regularity]
Let $\phi^0\in W$ be given. Then the solution to (\ref{e11})--(\ref{e12}) is regular in the Gevrey sense, see definition in Section 8, globally (locally resp.) in 2D (3D resp.), i.e., the solution is globally (locally resp.) spatially analytic in 2D (3D resp.).
\end{theorem}

For long time behavior of solutions to (\ref{e11})--(\ref{e12}), we have the following three theorems.

\begin{theorem}[2D Long time behavior for small initial perturbations]
Let $\phi^0\in V$ be given and assume that $\displaystyle{\bar{\phi}=\frac{1}{|\Omega|}\int_\Omega \phi^0(\mathbf{x})d\mathbf{x}>\sqrt{3}/3}$. Then there exists a constant $\delta>0$, sufficiently small, such that if $\|\phi^0-\bar{\phi}\|_{H^2}^2<\delta$ then the global solution to (\ref{e11})--(\ref{e12}) satisfies
$$
\|\phi(\cdot,t)-\bar{\phi}\|_{H^1}^2\le Ce^{-Ct},\ \ \forall\ t\ge0,
$$
for some constant $C>0$ independent of time.
\end{theorem}

\begin{theorem}[3D Long time behavior for small initial perturbations]
Let the condition of Theorem 1.6 be satisfied. Then the strong solution to (\ref{e11})--(\ref{e12}) exists globally in time and it satisfies
$$
\|\phi(\cdot,t)-\bar{\phi}\|_{H^1}^2\le Ce^{-Ct},\ \ \forall\ t\ge0,
$$
for some constant $C>0$ independent of time.
\end{theorem}

\begin{theorem}[2D Long time behavior for large initial perturbations]
Let $\phi^0\in V$ be given and suppose that $l\le L<\pi$. Then the global solution to (\ref{e11})--(\ref{e12}) satisfies
$$
\|\phi(\cdot,t)-\bar{\phi}\|_{H^2}^2\le Ce^{-Ct},\ \ \forall\ t\ge0,
$$
for some constant $C>0$ independent of time.
\end{theorem}

\begin{remark}
Our results show that the diffuse interface model is globally (locally resp.) well-posed in 2D (3D resp.). Moreover, under certain conditions on the initial data, the solution collapses to a constant equilibrium state as time tends to infinity. This suggests that the distinction between the tumor and the surrounding tissue in the microenvironment will blur as time proceeds. Theorems 1.6 and 1.7 indicate that constant solutions are locally asymptotically stable when they are outside the chemical spinodal. While, based on our analysis, it is unknown what happens if they are in the spinodal region where $F''(\bar{\phi})<0$. This launches a new interesting problem for future pursue. In addition, we expect to extend the results to the case with $S\neq0$ by adopting the ideas in this paper.
\end{remark}

We prove the above Theorem 1.1--Theorem 1.5 by combining the standard Galerkin
approximation and method of energy estimate.
The energy estimate is delicate mainly due
to the coupling between the Cahn-Hilliard equation and the
Hele-Shaw cell equations by advection and order parameter $\phi$. Great efforts have been made to
simplify the proofs. The proofs of Theorem 1.1--Theorem 1.5 involve intensive applications of
Sobolev and Gagliardo-Nirenberg type inequalities, see Lemma 2.2. The proofs of Theorem 1.6 and Theorem 1.7
are in the spirit of \cite{Boyer}. Since $\bar{\phi}$ is assumed to be outside the spinodal region ($\bar{\phi}>\sqrt{3}/3$), we note that $F''(\phi)=3\phi^2-1\ge 0$
(c.f. $(\ref{e11})_2$) on a small neighborhood $I_{\bar{\phi}}=[\bar{\phi}-\delta,\bar{\phi}+\delta]$, where $\delta>0$ is sufficiently small. The idea is to solve a modified problem with $F$ replaced by
an auxiliary function $F_{\bar{\phi}}$ whose second order derivative is non-negative on $\mathbb{R}$ and coincides with $F''(\cdot)$ on $I_{\bar{\phi}}$. Then, under the smallness assumption on the initial perturbations, it can be shown that the solution to the modified problem is indeed the solution to the original problem and converges exponentially to the constant state, as time goes to infinity. The last theorem is proved by
careful exploration of the condition $l\le L<\pi$ and previously established $a\ priori$ estimates. The proof involves exhaustive coupling of energy estimates.

The plan of the rest of this paper is as follows. In Section 2, we
give some basic facts that will be used in the proofs of Theorems 1.1--1.8. In Section 3, we present some lower order
$a\ priori$ estimates for the approximate solutions constructed through the Galerkin method. Then we prove Theorem 1.1 and
Theorem 1.2 in Section 4 and Section 5 respectively. In Section 6, we show the uniqueness of the solution by defining a stored
energy. Then we improve the spatial regularity of the solution and give a proof of Theorem 1.4 in Section 6. Section 7 is
devoted to the proof of the Gevrey regularity of the solution following ideas from \cite{FT,LT,OT}. Finally, we show the long time asymptotic behavior of the solution in Section 9.

\section{Preliminaries and functional setting}

In this section, we give some basic functional settings and preliminaries for the problem.
First we define the function space in which the approximate Galerkin solutions are established. For more details, we refer the readers to \cite{LT,OT} and references therein. For simplicity, we only give the 3D case, while the 2D case follows similarly. We let
$$
H_m=\mathrm{span}\{\phi_{i,j,k}|\lambda_{i,j,k}\le \Lambda_m\},
$$
where
\begin{equation}
\phi_{i,j,k}(\vec{\mathbf{x}})=\sqrt{\frac{8}{lLh}}\cos\Big(\frac{i\pi}{L}x\Big)\cos\Big(\frac{j\pi}{l}y\Big)\cos\Big(\frac{k\pi}{h}z\Big),\label{e21}
\end{equation}
for $i,j,k=0,1,2,...$ are the eigenfunctions of the Laplacian operator on the domain $\Omega=[0,L]\times[0,l]\times[0,h]$, subject to the boundary condition $\nabla\phi\cdot \mathbf{n}=0$ on $\partial \Omega$ with the corresponding eigenvalues
$$
\lambda_{i,j,k}=\Big(\frac{i\pi}{L}\Big)^2+\Big(\frac{j\pi}{l}\Big)^2+\Big(\frac{k\pi}{h}\Big)^2,
$$
and $0<\Lambda_1<\Lambda_2<\cdots<\Lambda_m<\cdots$ denote the set of distinct eigenvalues $\lambda_{i,j,k}$'s ordered by their magnitude.

Next, we define
$$
\widetilde{\mathcal{V}}:=\{\vec{\mathbf{v}}\in (C^\infty(\overline{\Omega}))^3|\vec{\mathbf{v}}\cdot\vec{\mathbf{n}}|_{\partial\Omega}=0\ \mathrm{and}\ \nabla\cdot\vec{\mathbf{v}}=0\ \mathrm{on}\ \Omega\},
$$
and let
$$
\mathrm{\mathbf{H}}:=\mathrm{Closure\ of\ }\widetilde{\mathcal{V}}\ \mathrm{in\ the}\ L^2\ -\ \mathrm{norm},
$$
$$
\mathrm{\mathbf {V}}:=\mathrm{Closure\ of\ }\widetilde{\mathcal{V}}\ \mathrm{in\ the}\ H^1\ -\ \mathrm{norm}.
$$

In addition, let
\begin{equation}
\vec{\mathbf{v}}_{i,j,k}=\left[\begin{aligned}
\frac{ik}{L}\sin\Big(\frac{i\pi}{L}x\Big)&\cos\Big(\frac{j\pi}{l}y\Big)\cos\Big(\frac{k\pi}{h}z\Big)\\
\frac{jk}{l}\cos\Big(\frac{i\pi}{L}x\Big)&\sin\Big(\frac{j\pi}{l}y\Big)\cos\Big(\frac{k\pi}{h}z\Big)\\
-\Big(\frac{i^2}{L^2}+\frac{j^2}{l^2}\Big)\cos\Big(&\frac{i\pi}{L}x\Big)\cos\Big(\frac{j\pi}{l}y\Big)\sin\Big(\frac{k\pi}{h}z\Big)
\end{aligned}\right],\label{e22}
\end{equation}
for $i,j=0,1,2,...;k=1,2,...$ and
\begin{equation}
\vec{\mathbf{w}}_{i,j,k}=\left[\begin{aligned}
\frac{j}{l}\sin\Big(\frac{i\pi}{L}x\Big)\cos&\Big(\frac{j\pi}{l}y\Big)\cos\Big(\frac{k\pi}{h}z\Big)\\
\frac{-i}{L}\cos\Big(\frac{i\pi}{L}x\Big)\sin&\Big(\frac{j\pi}{l}y\Big)\cos\Big(\frac{k\pi}{h}z\Big)\\
&0
\end{aligned}\right],\label{e23}
\end{equation}
for $i,j,k=0,1,2,...$. Then $\{\vec{\mathbf{v}}_{i,j,k}$'s,$\vec{\mathbf{w}}_{i,j,k}$'s$\}$ form an orthogonal basis for $\mathrm{\mathbf{H}}$. Moreover, the finite-dimensional subspace $\mathrm{\mathbf{H}}_m$ is defined in the similar fashion as $H_m$, using the orthogonal basis $\{\vec{\mathbf{v}}_{i,j,k}$'s,$\vec{\mathbf{w}}_{i,j,k}$'s$\}$ (c.f. \cite{F,LT,OT}).

Next, we let $\mathrm{\mathbf{P}}_\sigma:(L^2(\Omega))^3\rightarrow \mathrm{\mathbf{H}}$ be the Helmholtz-Leray orthogonal projection. By applying $\mathrm{\mathbf{P}}_\sigma$ to $(\ref{e11})_3$ we have
\begin{equation}\label{e24}
\vec{\mathbf{v}}-\mathrm{\mathbf{P}}_\sigma(\mu\nabla\phi)=0.
\end{equation}
The following lemma is given in \cite{LT} and references therein.

\begin{lemma}
For given functions $\mu$ and $\phi$, there exists a unique solution to the problem $(\ref{e11})_3$, $(\ref{e11})_4$ subject to the boundary condition $\vec{\mathbf{v}}\cdot\mathbf{n}|_{\partial\Omega}=0$. Moreover, there exists a constant $C>0$ such that
$$
\|\vec{\mathbf{v}}\|^2_{W^{s,q}}\le C\|\mu\nabla\phi\|_{W^{s,q}},\ \ s=0,1,\ \ \forall\ 1\le q<\infty.
$$
\end{lemma}

We also need the following Sobolev, Gagliardo-Nirenberg, Agmon and Poincar\'{e} type inequalities and interpolation inequalities which are standard and classic, c.f. \cite{Adam,CF,LSU,T}.

\begin{lemma} Let $\bar{f}=\int_\Omega f\ d\mathbf{x}$. Then in the 3D case, there exists a constant $C>0$ such that
$$
\begin{aligned}
(\mathrm{I})\ &\|f\|^2_{L^p}\le C\|f\|^2_{H^1},\ \ \forall\ p\in[1,6],\\
(\mathrm{II})\ &\|f\|^2_{L^4}\le C\big(\|f\|^{1/2}\|\nabla f\|^{3/2}+\|f\|^2\big),\\
(\mathrm{III})\ &\|f\|^2_{L^\infty}\le C\|f\|_{H^2}^2,\\
(\mathrm{IV})\ &\|f\|^2_{L^\infty}\le C\|f\|_{H^1}\|f\|_{H^2},\\
(\mathrm{V})\ &\|f\|^2_{L^2}\le C\|\nabla f\|_{L^2}^2,\ \ \mathrm{if}\ \ \bar{f}=0.
\end{aligned}
$$
In the 2D case, there exists a constant $C>0$ such that
$$
\begin{aligned}
(\mathrm{VI})\ &\|f\|^2_{L^4}\le C\big(\|f\|_{L^2}\|\nabla f\|_{L^2}+\|f\|_{L^2}^2\big),\\
(\mathrm{VII})\ &\|f\|^2_{L^4}\le C\|f\|_{L^2}\|\nabla f\|_{L^2},\ \ \mathrm{if}\ \ \bar{f}=0,\\
(\mathrm{VIII})\ &\|f\|^2_{L^\infty}\le C\|f\|_{L^2}\| f\|_{H^2}.
\end{aligned}
$$
In the 3D case, there exists a constant $C>0$ such that
$$
\begin{aligned}
(\mathrm{IX})\ &\|f\|^2_{L^3}\le C\big(\|f\|_{L^2}\|\nabla f\|_{L^2}+\|f\|^2_{L^2}\big),\\
(\mathrm{X})\ &\|f\|^2_{L^3}\le C\|f\|_{L^2}\|\nabla f\|_{L^2},\ \ \mathrm{if}\ \ \bar{f}=0,\\
(\mathrm{XI})\ &\|f\|^2_{L^\infty}\le C\|f\|_{L^2}\| f\|_{H^3},\\
(\mathrm{XII})\ &\|\nabla f\|_{L^2}^2\le C\|\Delta f\|^2_{L^2},\ \ \mathrm{if}\ \ \nabla f\cdot\mathbf{n}|_{\partial\Omega}=0,\\
(\mathrm{XIII})\ &\|f\|_{H^2}^2\le C\|\Delta f\|^2_{L^2},\ \ \mathrm{if}\ \ \bar{f}=0\ \ \mathrm{and}\ \ \nabla f\cdot\mathbf{n}|_{\partial\Omega}=0.
\end{aligned}
$$
\end{lemma}

Now we give the definition of strong solutions to the problem (\ref{e11})--(\ref{e12}).

\begin{definition}
Let $\phi^0\in V$ be given, and let $T>0$. A strong solution of (\ref{e11})--(\ref{e12}) in the
interval $[0,T]$ is a pair of functions $(\phi,\mu)$ such that $\phi\in C([0,T];V)\cap L^2([0,T];W)$, $\mu\in C([0,T];H)\cap L^2([0,T];V)$, $\displaystyle{\frac{d\phi}{dt}}\in L^2([0,T];H)$, and
$$
\frac{\partial\phi}{\partial t}+\vec{\mathbf{v}}\cdot\nabla\phi=\Delta \mu,\ \ \mathrm{in}\ H,
$$
for almost every $t\in[0,T]$, with $\mu$ given by $(\ref{e11})_2$ and $\vec{\mathbf{v}}$ given by (\ref{e24}).
That is, for every $\omega\in H$ the above equation holds in the following sense
\begin{equation}
\big(\phi(t_2),\omega\big)-\big(\phi(t_1),\omega\big)+\int^{t_2}_{t_1}\big(\vec{\mathbf{v}}\cdot\nabla\phi(\tau),\omega\big)d\tau=\int_{t_1}^{t_2}\big(\Delta \mu(\tau),\omega\big)d\tau,
\end{equation}\label{ss}
for every $t_1,t_2\in[0,T]$, where $\mu$ is given by $(\ref{e11})_2$ and $\vec{\mathbf{v}}$ is given by (\ref{e24}).
\end{definition}

\section{Lower order estimates in 2D and 3D}

In this section, we will give some general energy estimates which are valid in both 2D and 3D. The energy estimates are
performed for the approximate solutions to the original problem, which are constructed through the standard Galerkin procedure (c.f. \cite{LT,OT}).

We consider the Galerkin approximation system:
\begin{equation}
\frac{d\phi_m}{dt}+P_m(\vec{\mathbf{v}}_m\cdot\nabla\phi_m)=\Delta\mu_m,\label{e31}
\end{equation}
\begin{equation}
\mu_m=P_m(\phi_m^3)-\phi_m-\Delta\phi_m,\label{e32}
\end{equation}
\begin{equation}
\vec{\mathbf{v}}_m=\mathrm{\mathbf{P}}_m\mathrm{\mathbf{P}}_\sigma(\mu_m\nabla\phi_m),\label{e33}
\end{equation}
with initial data
\begin{equation}
\phi_m(\vec{\mathbf{x}},0)=\phi_m^0(\vec{\mathbf{x}})=P_m\phi^0(\vec{\mathbf{x}}),\label{e34}
\end{equation}
where $\vec{\mathbf{v}}_m\in \mathrm{\mathbf{H}}_m$, $\phi_m,\mu_m\in H_m$, and $\mathrm{\mathbf{P}}_m$ and $P_m$ denote the orthogonal projections from $\mathrm{\mathbf{H}}$ and $H$ onto $\mathrm{\mathbf{H}}_m$ and $H_m$, respectively. The finite-dimensional ODE system has a short time unique solution due to its locally Lipschitz nonlinearity.

We now prove some lower order energy estimates which are valid in 2D and 3D.

{\bf Step 1.}
Taking the $L^2$ inner product of (\ref{e31}) with $\mu_m$ we have
\begin{equation}
\frac{d}{dt}\bigg(\int_\Omega\Big(\frac{1}{4}\phi_m^4-\frac{1}{2}\phi_m^2\Big)d\mathbf{x}+\frac{1}{2}\|\nabla\phi_m\|^2\bigg)+\|\nabla\mu_m\|^2
=-\big(P_m(\vec{\mathbf{v}}_m\cdot\nabla\phi_m),\mu_m\big).\label{e35}
\end{equation}
By definition, it is straightforward to show that
\begin{equation}
-\big(P_m(\vec{\mathbf{v}}_m\cdot\nabla\phi_m),\mu_m\big)=-\big(\vec{\mathbf{v}}_m\cdot\nabla\phi_m,\mu_m \big),\label{e36}
\end{equation}
which, together with (\ref{e35}), implies that
\begin{equation}
\frac{d}{dt}\bigg(\int_\Omega\Big(\frac{1}{4}\phi_m^4-\frac{1}{2}\phi_m^2\Big)d\mathbf{x}+\frac{1}{2}\|\nabla\phi_m\|^2\bigg)+\|\nabla\mu_m\|^2
=-\big(\vec{\mathbf{v}}_m\cdot\nabla\phi_m,\mu_m \big).\label{e37}
\end{equation}

Next, taking the $L^2$ inner product of (\ref{e33}) with $\vec{\mathbf{v}}_m$ we have
\begin{equation}
\|\vec{\mathbf{v}}_m\|^2=\big(\mathrm{\mathbf{P}}_m\mathrm{\mathbf{P}}_\sigma(\mu_m\nabla\phi_m), \vec{\mathbf{v}}_m \big).\label{e38}
\end{equation}
Again, by definition we have
\begin{equation}
\big(\mathrm{\mathbf{P}}_m\mathrm{\mathbf{P}}_\sigma(\mu_m\nabla\phi_m), \vec{\mathbf{v}}_m \big)=\big(\mu_m\nabla\phi_m,\vec{\mathbf{v}}_m \big)=
\big(\vec{\mathbf{v}}_m\cdot\nabla\phi_m,\mu_m \big),\label{e39}
\end{equation}
which yields
\begin{equation}
\|\vec{\mathbf{v}}_m\|^2=\big(\vec{\mathbf{v}}_m\cdot\nabla\phi_m,\mu_m \big).\label{e310}
\end{equation}

Therefore, by adding (\ref{e37}) and (\ref{e310}) we have
\begin{equation}
\frac{d}{dt}\bigg(\int_\Omega\Big(\frac{1}{4}\phi_m^4-\frac{1}{2}\phi_m^2\Big)d\mathbf{x}+\frac{1}{2}\|\nabla\phi_m\|^2\bigg)+\|\nabla\mu_m\|^2+\|\vec{\mathbf{v}}_m\|^2
=0,\label{e311}
\end{equation}
which is equivalent to
\begin{equation}
\frac{d}{dt}\bigg(\frac{1}{4}\|\phi_m^2-1\|^2+\frac{1}{2}\|\nabla\phi_m\|^2\bigg)+\|\nabla\mu_m\|^2+\|\vec{\mathbf{v}}_m\|^2
=0,\label{e312}
\end{equation}
and this in turns implies, by integrating w.r.t. time, that
\begin{equation}
\begin{aligned}
\bigg(\frac{1}{4}\|\phi_m^2-1\|^2+\frac{1}{2}\|\nabla\phi_m\|^2\bigg)(t)&+\int_0^t\big(\|\nabla\mu_m\|^2+\|\vec{\mathbf{v}}_m\|^2\big)ds\\
=&\bigg(\frac{1}{4}\|\phi_m^2(0)-1\|^2+\frac{1}{2}\|\nabla\phi_m(0)\|^2\bigg)\le C(\Omega,\|\phi^0\|_{H^1}),
\end{aligned}\label{e313}
\end{equation}
where $C(\Omega,\|\phi^0\|_{H^1})$ denotes a constant which depends on $\Omega$ and $\|\phi^0\|_{H^1}$. In the last step we used the fact that $\|\phi_m(0)\|_{L^4}\le C\|\phi_m(0)\|_{H^1}\le C\|\phi^0\|_{H^1}$ (see Lemma 2.2).
As a result of (\ref{e313}), it follows that (\ref{e31})--(\ref{e34}) exists globally in time, for positive time, because the nonlinearity
is locally Lipschitz and the solution is bounded.

{\bf Step 2.} In this step, we derive some by-products of (\ref{e313}). Integrating equation (\ref{e31}) over $\Omega\times [0,t]$ we have
\begin{equation}
\int_\Omega\phi_m(\mathbf{x},t)d\mathbf{x}=\int_\Omega\phi_m^0(\mathbf{x})d\mathbf{x}=\int_\Omega\phi^0(\mathbf{x})d\mathbf{x}.\label{e314}
\end{equation}
Define
\begin{equation}
\bar{\phi}_m=\bar{\phi}\equiv\frac{1}{|\Omega|}\int_\Omega\phi^0(\mathbf{x})d\mathbf{x}.\label{e315}
\end{equation}
From (\ref{e313}) we easily see that
\begin{equation}
\|\phi_m\|_{H^1}^2=\|\phi_m\|^2+\|\nabla\phi_m\|^2\le |\Omega|^{1/2}\|\phi_m\|_{L^4}^2+\|\nabla\phi_m\|^2\le C\|\phi_m^2-1\|^2+\|\nabla\phi_m\|^2\le C(\Omega,\|\phi^0\|_{H^1}),\label{e316}
\end{equation}
which implies that
\begin{equation}
\|\phi_m\|_{H^1}^3\le C(\Omega,\|\phi^0\|_{H^1}).\label{e317}
\end{equation}
Then we have, by virtue of Lemma 2.2 (I)
\begin{equation}
\bigg|\int_\Omega\phi^3_m(\mathbf{x},t)d\mathbf{x}\bigg|\le |\Omega|^{1/2}\|\phi_m\|_{L^6}^3\le C\|\phi_m\|_{H^1}^3\le C(\Omega,\|\phi^0\|_{H^1})\le C.\label{e318}
\end{equation}

Thanks to the boundary condition $\nabla\phi_m\cdot\mathbf{n}|_{\partial\Omega}=0$ and (\ref{e314})
\begin{equation}
\bar{\mu}_m\equiv\frac{1}{|\Omega|}\int_\Omega\mu_m(\mathbf{x},t)d\mathbf{x}=\frac{1}{|\Omega|}\int_\Omega\big(
\phi^3_m(\mathbf{x},t)-\phi_m-\Delta\phi_m\big)d\mathbf{x}=\frac{1}{|\Omega|}\int_\Omega
\phi^3_m(\mathbf{x},t)d\mathbf{x}-\bar{\phi}.\label{e319}
\end{equation}
Then, from (\ref{e318}) we know that for any $t\ge0$ and $m\ge1$, it holds that
\begin{equation}
|\bar{\mu}_m|\le C.\label{e320}
\end{equation}

Due to the boundary conditions and the trigonometric basis (\ref{e21})--(\ref{e23}), by using Poincar\'{e}'s inequality one can easily show that
\begin{equation}
\|\phi_m-\bar{\phi}_m\|_{H^{2s}}\le C\|(-\Delta)^s\phi_m\|,\ \ \|\phi_m-\bar{\phi}_m\|_{H^{2s+1}}\le C\|\nabla(-\Delta)^s\phi_m\|,\ \ s\ge1,\label{e321}
\end{equation}
and
\begin{equation}
\|\mu_m-\bar{\mu}_m\|_{H^{2s}}\le C\|(-\Delta)^s\mu_m\|,\ \ \|\mu_m-\bar{\mu}_m\|_{H^{2s+1}}\le C\|\nabla(-\Delta)^s\mu_m\|,\ \ s\ge1.\label{e322}
\end{equation}

Furthermore, if we write $\hat{\phi}_m=\phi_m-\bar{\phi}$ and $\hat{\mu}_m=\mu_m-\bar{\mu}_m$, then (\ref{e31})--(\ref{e33}) is equivalent to
\begin{equation}\label{ns}
\left\{
\begin{aligned}
&\frac{d\hat{\phi}_m}{dt}+P_m(\vec{\mathbf{v}}_m\cdot\nabla\hat{\phi}_m)=\Delta\hat{\mu}_m,\\
&\hat{\mu}_m=P_m\big((\hat{\phi}_m+\bar{\phi})^3\big)-\hat{\phi}_m-\bar{\phi}-\Delta\hat{\phi}_m-\bar{\mu}_m\\
&\ \ \ \ =P_m\big[(\hat{\phi}_m+\bar{\phi})^3-\overline{(\hat{\phi}_m+\bar{\phi})^3}\big]-\hat{\phi}_m-\Delta\hat{\phi}_m,\\
&\vec{\mathbf{v}}_m=\mathrm{\mathbf{P}}_m\mathrm{\mathbf{P}}_\sigma(\hat{\mu}_m\nabla\hat{\phi}_m).
\end{aligned}
\right.
\end{equation}
Notice that in order to obtain the last equation we used the fact that $\bar{\mu}_m$ is $\mathbf{x}$-independent and that $\mathrm{\mathbf{P}}_\sigma(\nabla{\phi}_m)=0$. From now on we will call $\hat{\phi}_m$ and $\hat{\mu}_m$ by $\phi_m$ and $\mu_m$, respectively, and use the above system for the next estimates, taking into considerations that $\overline{\hat{\phi}_m}=0$, $\overline{\hat{\mu}_m}=0$, and using the already established estimates (\ref{e316}) and (\ref{e320}).

{\bf Step 3.} Taking the $L^2$ inner product of $(\ref{ns})_1$ with $\phi_m$ and using (\ref{e313}) we have
\begin{equation}
\frac{1}{2}\frac{d}{dt}\|\phi_m\|^2+3\|(\phi_m+\bar{\phi})\nabla\phi_m\|^2+\|\Delta\phi_m\|^2=\|\nabla\phi_m\|^2\le C(\Omega,\phi^0),\label{e323}
\end{equation}
where $C(\Omega,\phi^0)$ denotes a constant depending only on $\Omega$ and initial data $\phi^0$, and we have used the fact that $(P_m(\vec{\mathbf{v}}_m\cdot\nabla\phi_m),\phi_m)=0$. Integrating (\ref{e323}) in time we have
\begin{equation}
\begin{aligned}
\frac{1}{2}\sup_{\tau\in[0,t]}\|\phi_m(\tau)\|^2+\int_0^t\big(3\|(\phi_m+\bar{\phi})\nabla\phi_m\|^2+\|\Delta\phi_m\|^2\big)ds &\le C(\Omega,\phi^0)t+\frac12\|\phi_m(0)\|^2\\
&\le C(\Omega,\phi^0)(t+1).\label{e324}
\end{aligned}
\end{equation}
Especially, by (\ref{e321}) and (\ref{e324}) we have
\begin{equation}
\int_0^t\|\phi_m\|^2_{H^2}ds\le C.\label{e325}
\end{equation}

{\bf Step 4.} As consequences of (\ref{e313}) and (\ref{e325}) we have the following. First, by $(\ref{ns})_2$ we have
\begin{equation}
\begin{aligned}
\|\nabla\Delta\phi_m\|^2&=\|\nabla P_m((\phi_m+\bar{\phi})^3)-\nabla\phi_m-\nabla\mu_m\|^2\\
&\le C\big(\|\nabla((\phi_m+\bar{\phi})^3)\|^2+\|\nabla\phi_m\|^2+\|\nabla\mu_m\|^2\big)\\
&\le C\big(\|(\phi_m+\bar{\phi})^2\nabla\phi_m\|^2+\|\nabla\phi_m\|^2+\|\nabla\mu_m\|^2\big)\\
&\le C\big(\|\phi_m^2\nabla\phi_m\|^2+\|\nabla\phi_m\|^2+\|\nabla\mu_m\|^2\big).
\end{aligned}\label{e326}
\end{equation}
By H\"{o}lder's inequality, Lemma 2.2 (I) and (\ref{e321}) we then have
\begin{equation}
\begin{aligned}
\mathrm{R.H.S.\ of}\ (\ref{e326})&\le C\big(\|\phi_m\|_{L^\infty}^4\|\nabla\phi_m\|^2+\|\nabla\phi_m\|^2+\|\nabla\mu_m\|^2\big)\\
&\le C\big(\|\nabla\phi_m\|^2\|\Delta\phi_m\|^2\|\nabla\phi_m\|^2+\|\nabla\phi_m\|^2+\|\nabla\mu_m\|^2\big)\\
&\le C\big(\|\nabla\phi_m\|^4\|\Delta\phi_m\|^2+\|\nabla\phi_m\|^2+\|\nabla\mu_m\|^2\big)\\
&\le C\big(\|\Delta\phi_m\|^2+\|\nabla\mu_m\|^2\big),
\end{aligned}\label{e327}
\end{equation}
where we used Lemma 2.2 (IV) and also applied (\ref{e313}) for the estimate of $\|\nabla\phi_m\|^4$. By plugging (\ref{e327}) into (\ref{e326}) we have
\begin{equation}
\|\nabla\Delta\phi_m\|^2\le C\big(\|\Delta\phi_m\|^2+\|\nabla\mu_m\|^2\big),\label{e328}
\end{equation}
which, together with (\ref{e313}) and (\ref{e325}), yields
\begin{equation}
\int_0^t\|\phi_m\|^2_{H^3}ds\le C.\label{e329}
\end{equation}
We note that the constant $C$ on the RHS of (\ref{e328}) is independent of time.

Second, similar to (\ref{e326})--(\ref{e327}) we have
\begin{equation}
\begin{aligned}
\|\mu_m\|^2&\le C\big(\|(\phi_m+\bar{\phi})^3-\overline{(\phi_m+\bar{\phi})^3}\|^2+\|\phi_m\|^2+\|\Delta\phi_m\|^2\big)\\
&\le C\big(\|\nabla(\phi_m+\bar{\phi})^3\|^2+\|\nabla\phi_m\|^2+\|\Delta\phi_m\|^2\big)\\
&\le C\big(\|(\phi_m+\bar{\phi})^2\nabla\phi_m\|^2+\|\Delta\phi_m\|^2\big)\\
&\le C\|\Delta\phi_m\|^2,
\end{aligned}\label{e43}
\end{equation}
where we used Lemma 2.2 (XII). It also holds that
\begin{equation}
\begin{aligned}
\|\nabla\mu_m\|^2&\le C\big(\|(\phi_m+\bar{\phi})^2\nabla\phi_m\|^2+\|\nabla\phi_m\|^2+\|\nabla\Delta\phi_m\|^2\big)\\
&\le C\big(\|\Delta\phi_m\|^2+\|\nabla\Delta\phi_m\|^2\big)\\
&\le C\|\nabla\Delta\phi_m\|^2.
\end{aligned}\label{e44}
\end{equation}

We remark that, according to the definition of strong solutions, with the estimates obtained above, it is still not enough to take limit, in $m$, and conclude the existence of solutions to the original problem due to the nonlinearity in the equations. Therefore, we have to seek higher order estimate of the approximate solutions.

{\bf Step 5.} Taking the $L^2$ inner product of $(\ref{ns})_1$ with $\Delta^2\phi_m$ and applying the Cauchy-Schwarz and Young inequalities we have
\begin{equation}
\begin{aligned}
&\frac{1}{2}\frac{d}{dt}\|\Delta\phi_m\|^2+\|\Delta^2\phi_m\|^2\\=&-\int_\Omega (\vec{\mathbf{v}}_m\cdot\nabla\phi_m)(\Delta^2\phi_m) d\mathbf{x}+
\int_\Omega\Delta((\phi_m+\bar{\phi})^3)(\Delta^2\phi_m)d\mathbf{x}-\int_\Omega\Delta\phi_m(\Delta^2\phi_m)d\mathbf{x}\\
\le &\frac{1}{2}\|\Delta^2\phi_m\|^2+\frac{3}{2}\|(\vec{\mathbf{v}}_m\cdot\nabla\phi_m)\|^2+\frac{3}{2}\|\Delta((\phi_m+\bar{\phi}^3))\|^2+\frac{3}{2}\|\Delta\phi_m\|^2.
\end{aligned}\label{e330}
\end{equation}
After rearranging terms we update above estimate as
\begin{equation}
\begin{aligned}
\frac{d}{dt}\|\Delta\phi_m\|^2+\|\Delta^2\phi_m\|^2
\le 3\|(\vec{\mathbf{v}}_m\cdot\nabla\phi_m)\|^2+3\|\Delta((\phi_m+\bar{\phi})^3)\|^2+3\|\Delta\phi_m\|^2.
\end{aligned}\label{e331}
\end{equation}

For the second term on the RHS of (\ref{e331}), by a direct calculation we have
\begin{equation}
\begin{aligned}
&\|\Delta((\phi_m+\bar{\phi})^3)\|^2\\
\le &C\big(\|(\phi_m+\bar{\phi})\|_{L^\infty}^4\|\Delta\phi_m\|^2+\|(\phi_m+\bar{\phi})\|^2_{L^\infty}\|\nabla\phi_m\|_{L^4}^4\big)\\
\le &C\big(\|\phi_m\|_{L^\infty}^4+1\big)\|\Delta\phi_m\|^2+C\big(\|\phi_m\|^2_{L^\infty}+1\big)\|\nabla\phi_m\|_{L^4}^4.
\end{aligned}\label{e332}
\end{equation}
Using Lemma 2.2 (IV), (\ref{e313}) and (\ref{e321}) we estimate the first term on the RHS of (\ref{e332}) as
\begin{equation}
\begin{aligned}
&C\big(\|\phi_m\|_{L^\infty}^4+1\big)\|\Delta\phi_m\|^2\\
\le &C\big(\|\phi_m\|_{H^1}^2\|\phi_m\|_{H^2}^2+1\big)\|\Delta\phi_m\|^2\\
\le &C\big(\|\Delta\phi_m\|^2+1\big)\|\Delta\phi_m\|^2\\
\le &C\big(\|\Delta\phi_m\|^4+\|\Delta\phi_m\|^2\big).
\end{aligned}\label{e333}
\end{equation}
Similarly, by using Lemma 2.2 (II) we have
\begin{equation}
\begin{aligned}
&C\big(\|\phi_m\|^2_{L^\infty}+1\big)\|\nabla\phi_m\|_{L^4}^4\\
\le &C\big(\|\phi_m\|_{H^1}\|\phi_m\|_{H^2}+1\big)\big(\|\nabla\phi_m\|\|\Delta\phi_m\|^3+\|\nabla\phi_m\|^4\big)\\
\le &C\big(\|\Delta\phi_m\|+1\big)\big(\|\Delta\phi_m\|^3+\|\nabla\phi_m\|^3\big)\\
\le &C\big(\|\Delta\phi_m\|+1\big)\|\Delta\phi_m\|^3\\
\le &C\big(\|\Delta\phi_m\|^4+\|\Delta\phi_m\|^2\big).
\end{aligned}\label{e334}
\end{equation}
Plugging (\ref{e333}) and (\ref{e334}) into (\ref{e332}) we have
\begin{equation}
\|\Delta((\phi_m+\bar{\phi})^3)\|^2\le C\big(\|\Delta\phi_m\|^4+\|\Delta\phi_m\|^2\big).\label{e335}
\end{equation}
So we update (\ref{e331}) as
\begin{equation}
\frac{d}{dt}\|\Delta\phi_m\|^2+\|\Delta^2\phi_m\|^2
\le 3\|(\vec{\mathbf{v}}_m\cdot\nabla\phi_m)\|^2+C\big(\|\Delta\phi_m\|^4+\|\Delta\phi_m\|^2\big).\label{e336}
\end{equation}

We remark that so far we have been focusing on the 3D case with all the estimates. Next, we shall deal with the first term on the RHS of (\ref{e336}). Due to differences between Sobolev embeddings in 2D and 3D, we have to separate the arguments into the two and three space dimensions, which will be given in the next two sections.

\section{Global existence in 2D}

In this section we prove Theorem 1.1, i.e., global existence of strong solutions to (\ref{e11})--(\ref{e12}) in 2D, utilizing of the $a\ priori$ energy estimates established in previous section.

First, by (\ref{e313}) and (\ref{e321}) we have
\begin{equation}
\begin{aligned}
\|(\vec{\mathbf{v}}_m\cdot\nabla\phi_m)\|^2&\le \|\vec{\mathbf{v}}_m\|^2\|\nabla\phi_m\|^2_{L^\infty}\\
&\le C\|\mu_m\nabla\phi_m\|^2\|\nabla\phi_m\|\|\nabla\phi_m\|_{H^2}\\
&\le C\|\mu_m\|^2_{L^4}\|\nabla\phi_m\|^2_{L^4}\|\nabla\Delta\phi_m\|.
\end{aligned}\label{e41}
\end{equation}
We estimate the RHS of (\ref{e41}) as follows. For the first term, using Lemma 2.2 (VII) and (\ref{e43})--(\ref{e44}) we have
\begin{equation}
\begin{aligned}
\|\mu_m\|^2_{L^4}\le C\|\mu_m\|\|\nabla\mu_m\|\le C\|\Delta\phi_m\|\|\nabla\Delta\phi_m\|.
\end{aligned}\label{e42}
\end{equation}
Since $\nabla\phi_m\cdot\mathbf{n}|_{\partial\Omega}=0$, the second term is estimated, using Lemma 2.2 (VI) and (XII) and (\ref{e313}), as
\begin{equation}
\begin{aligned}
\|\nabla\phi_m\|^2_{L^4}\le C\|\nabla\phi_m\|\|\Delta\phi_m\|\le C\|\Delta\phi_m\|.
\end{aligned}\label{e45}
\end{equation}

Combining (\ref{e41}), (\ref{e42}) and (\ref{e45}) we have
\begin{equation}
\begin{aligned}
\|(\vec{\mathbf{v}}_m\cdot\nabla\phi_m)\|^2\le C\|\Delta\phi_m\|^2\|\nabla\Delta\phi_m\|^2.
\end{aligned}\label{e411}
\end{equation}


Plugging (\ref{e411}) into (\ref{e336}) and applying Lemma 2.2 (V) to $\Delta\phi_m$ we have
\begin{equation}
\begin{aligned}
\frac{1}{2}\frac{d}{dt}\|\Delta\phi_m\|^2+\frac{1}{2}\|\Delta^2\phi_m\|^2
&\le C\|\nabla\Delta\phi_m\|^2\|\Delta\phi_m\|^2+C\big(\|\Delta\phi_m\|^4+\|\Delta\phi_m\|^2\big)\\
&\le C\|\nabla\Delta\phi_m\|^2\|\Delta\phi_m\|^2+C\|\nabla\Delta\phi_m\|^2.
\end{aligned}\label{e412}
\end{equation}
Applying Gronwall's inequality to (\ref{e412}) and using (\ref{e329}) we have
\begin{equation}
\sup_{\tau\in[0,t]}\|\phi_m(\tau)\|^2_{H^2}+\int_0^t\|\phi_m\|^2_{H^4}ds\le C.\label{e413}
\end{equation}

As consequences of (\ref{e413}) we have
\begin{equation}
\begin{aligned}
\|\mu_m\|^2_{H^2}&\le C\|\Delta\mu_m\|^2\\
&\le C\big(\|\Delta((\phi_m+\bar{\phi})^3)\|^2+\|\Delta\phi_m\|^2+\|\Delta^2\phi_m\|^2\big)\\
&\le C\big(\|\Delta\phi_m\|^4+\|\Delta\phi_m\|^2+\|\Delta^2\phi_m\|^2\big),\ \ \ \mathrm{by}\ (\ref{e335}),\\
&\le C\big(\|\Delta\phi_m\|^4+\|\Delta^2\phi_m\|^2\big),\ \ \ \mathrm{by}\ (\ref{e321}),
\end{aligned}\label{e414}
\end{equation}
which implies, by (\ref{e413}) that
\begin{equation}
\int_0^t\|\mu_m\|^2_{H^2}ds\le C.\label{e415}
\end{equation}

From (\ref{e411}) and (\ref{e413}) we see that
\begin{equation}
\int_0^t\|(\vec{\mathbf{v}}_m\cdot\nabla\phi_m)\|^2ds\le C\int_0^t \|\Delta\phi_m\|^2\|\nabla\Delta\phi_m\|^2ds\le C\int_0^t \|\nabla\Delta\phi_m\|^2ds\le C.\label{e416}
\end{equation}

Combining (\ref{e415})--({\ref{e416}}) and $(\ref{ns})_1$ we have
\begin{equation}
\int_0^t\Big\|\frac{d\phi_m}{dt}\Big\|^2ds\le C.\label{e417}
\end{equation}

Moreover, we have
\begin{equation}
\begin{aligned}
\int_0^t\|\vec{\mathbf{v}}_m\|_{H^1}^2ds &\le \int_0^t\big(\|\mu_m\nabla\phi_m\|^2+\|\nabla\mu_m\nabla\phi_m\|^2+\|\mu_m\nabla^2\phi_m\|^2\big)ds\\
&\le \int_0^t\big(\|\mu_m\|_{L^\infty}^2\|\nabla\phi_m\|^2+\|\nabla\mu_m\|_{L^4}^2\|\nabla\phi_m\|^2_{L^4}+\|\mu_m\|_{L^\infty}^2\|\Delta\phi_m\|^2\big)ds\\
&\le C\int_0^t\|\mu_m\|_{H^2}^2\|\Delta\phi_m\|^2ds\\
&\le C\int_0^t\|\mu_m\|_{H^2}^2ds\le C.
\end{aligned}\label{e418}
\end{equation}

Collecting above estimates we have
\begin{equation}
\sup_{\tau\in[0,t]}\big(\|\phi_m(\tau)\|^2_{H^2}+\|\mu_m(\tau)\|^2\big)+
\int_0^{t}\big(\|\vec{\mathbf{v}}_m\|^2_{H^1}+\|\phi_m\|^2_{H^4}+\|\mu_m\|^2_{H^2}+\Big\|\frac{d\phi_m}{dt}\Big\|^2\big)ds\le C(t),\ \ \forall\ t\ge0.\label{e419}
\end{equation}

Now, we are ready to use these estimates and apply standard compactness theorems to pass to the limit, in $m$, in order to show the global existence result in the 2D case. Fix $T>0$. In what follows, all arguments will be carried out within the interval $[0,T]$ unless specified otherwise. First, we note that, by (\ref{e419}), the sequence $\{\phi_m\}$ is bounded in $L^\infty([0,T];H^2)$ and $L^2([0,T];H^4)$ respectively, then by the weak compactness theorem, there exists a $\phi\in L^\infty([0,T];H^2)\cap L^2([0,T];H^4)$ and a subsequence $\{\phi_{m'}\}$ of $\{\phi_m\}$ such that $\{\phi_{m'}\}$ converges to $\phi$ in the weak* and weak topology of $L^\infty([0,T];H^2)$ and $L^2([0,T];H^4)$, respectively. Then we have from (\ref{e419}) that
\begin{equation}
\|\phi(s)\|_{H^2}\le \liminf_{m\to\infty}\|\phi_m(s)\|_{H^2}\le C,\ \ \mathrm{for\ a.e.}\ s\in[0,t],
\end{equation}
and
\begin{equation}
\int_0^t\|\phi(s)\|^2_{H^4}ds\le \liminf_{m\to\infty}\int_0^t\|\phi_m(s)\|^2_{H^4}ds\le C.
\end{equation}
In addition, since $H^4$ is compactly embedded in $H^3$, and since the sequence $\{\frac{d\phi_m}{dt}\}$ is bounded in $L^2([0,T];L^2)$, by Aubin's Compactness Theorem we know that $\{\phi_{m'}\}$ converges strongly to $\phi$ in $L^2([0,T];H^3)$.
From now on, we denote all the subsequences of $\{\phi_m\}$ by itself by extracting and relabeling. Thus, $\{\phi_m\}$ also converges strongly to $\phi$ in $H^3$ almost everywhere in $[0,T]$. In particular,
\begin{equation}\label{E}
\begin{aligned}
\|\phi_m(t)\|_{H^3}\longrightarrow \|\phi(t)\|_{H^3}\ &\ \mathrm{pointwise\ everywhere\ on}\ E\subseteq[0,T];\\
&\big|[0,T]\backslash E\big|=0.
\end{aligned}
\end{equation}
By virtue of the estimate (\ref{e419}) we know that there exists a function $\vec{\mathbf{v}}\in L^2([0,T];\mathrm{\mathbf{V}})$ such that $\{\vec{\mathbf{v}}_m\}$ converges weakly to $\vec{\mathbf{v}}$ in $L^2([0,T];\mathrm{\mathbf{V}})$. Therefore, we have
\begin{equation}
\int_0^T\|\vec{\mathbf{v}}(s)\|^2_{\mathrm{\mathbf{V}}}ds\le \liminf_{m\to\infty}\int_0^T\|\vec{\mathbf{v}}_m(s)\|^2_{\mathrm{\mathbf{V}}}ds\le C.
\end{equation}
In summary, we have
\begin{equation}\label{con}
\begin{aligned}
\phi_m &\longrightarrow \phi\ \ \mathrm{weakly\ in}\ \ L^2([0,T];H^4),\\
\phi_m &\longrightarrow \phi\ \ \mathrm{strongly\ in}\ \ L^2([0,T];H^3),\\
\phi_m &\longrightarrow \phi\ \ \mathrm{weak^*\ in}\ \ L^\infty([0,T];H^2),\\
\frac{d\phi_m}{dt} &\longrightarrow \frac{d\phi}{dt}\ \ \mathrm{weakly\ in}\ \ L^2([0,T];L^2),\\
\vec{\mathbf{v}}_m &\longrightarrow \vec{\mathbf{v}}\ \ \mathrm{weakly\ in}\ \ L^2([0,T];\mathrm{\mathbf{V}}),
\end{aligned}
\end{equation}
respectively.

Since the sequence $\{\phi_m\}$ is bounded in $L^\infty([0,T];H^2)$ and $\{\frac{d\phi_m}{dt}\}$ is bounded in $L^2([0,T];L^2)$, we know from the Arzela-Ascoli Theorem
that $\{\phi_m\}$ converges to $\phi$ in $C([0,T];L^2)$, and hence in $C_W([0,T];L^2)$, $C_W([0,T];H^1)$ and $C_W([0,T];H^2)$. 
In particular, we have
\begin{equation}\label{e419'}
(\phi_m(t_2),\omega)-(\phi_m(t_1),\omega)\longrightarrow (\phi(t_2),\omega)-(\phi(t_1),\omega),
\end{equation}
for every $t_1,t_2\in[0,T]$ and every $\omega\in H$. Moreover, from the first weak convergence in (\ref{con}) we know that for every $t_1,t_2\in[0,T]$ and every $\omega\in H$,
\begin{equation}\label{e419''}
\int_{t_1}^{t_2}(-\Delta^2\phi_m-\Delta\phi_m,\omega)ds\longrightarrow \int_{t_1}^{t_2}(-\Delta^2\phi-\Delta\phi,\omega)ds.
\end{equation}

Next, we will show that
\begin{equation}\label{e420}
\int_{t_1}^{t_2}(\Delta P_m(\phi_m^3),\omega)ds-\int_{t_1}^{t_2}(\Delta (\phi^3),\omega)ds \longrightarrow 0,
\end{equation}
as $m\to\infty$ for every $t_1,t_2\in[0,T]$ and every $\omega\in H$. First, we notice that, by Lemma 2.2 (III) and $(\ref{con})_3$
\begin{equation}\label{e421}
\begin{aligned}
\|\phi_m^3-\phi^3\|&=\|(\phi_m-\phi)(\phi_m^2+\phi_m\phi+\phi^2)\|,\\
&\le 2\big(\|\phi_m\|^2_{L^\infty}+\|\phi\|^2_{L^\infty}\big)\|\phi_m-\phi\|\\
&\le 2\big(\|\phi_m\|^2_{H^2}+\|\phi\|^2_{H^2}\big)\|\phi_m-\phi\|\\
&\le C\|\phi_m-\phi\|.
\end{aligned}
\end{equation}
Second, by a simple calculation and (\ref{e419}) one can show that
\begin{equation}\label{e422}
\begin{aligned}
\|\Delta^2(\phi_m^3)\|&\le C\big(\|\phi_m\|^2_{L^\infty}\|\Delta^2\phi_m\|+\|\phi_m\|_{L^\infty}\|\nabla\phi_m\|_{L^4}\|\nabla\Delta\phi_m\|_{L^4}\\
&\ \ \ \ \ \ +\|\nabla\phi_m\|_{L^6}^2\|\Delta\phi_m\|_{L^6}+\|\phi_m\|_{L^\infty}\|\Delta\phi_m\|_{L^\infty}\|\Delta\phi_m\|\big)\\
&\le C\big(\|\phi_m\|^2_{H^2}\|\Delta^2\phi_m\|+\|\phi_m\|_{H^2}\|\nabla\phi_m\|_{H^1}\|\nabla\Delta\phi_m\|_{H^1}\\
&\ \ \ \ \ \ +\|\nabla\phi_m\|_{H^1}^2\|\Delta\phi_m\|_{H^1}+\|\phi_m\|_{H^2}\|\Delta\phi_m\|_{H^2}\|\Delta\phi_m\|\big)\\
&\le C\|\Delta^2\phi_m\|.
\end{aligned}
\end{equation}
And the same is true for $\phi$, i.e., $\|\Delta^2(\phi^3)\|\le C\|\Delta^2\phi\|$. Then we have
\begin{equation}\label{e423}
\begin{aligned}
\|\Delta P_m(\phi_m^3)-\Delta (\phi^3)\|^2&=\Big(P_m(\phi_m^3)-\phi^3,\Delta^2\big(P_m(\phi_m^3)-\phi^3\big)\Big)\\
&\le \|P_m(\phi_m^3)-\phi^3\|_{L^2}\|\Delta^2\big(P_m(\phi_m^3)-\phi^3\big)\|_{L^2}\\
&\le C\big(\|\phi_m-\phi\|_{L^2}+\|(P_m-I)\phi^3\|_{L^2}\big)\big(\|\Delta^2\phi_m\|_{L^2}+\|\Delta^2\phi\|_{L^2}\big),
\end{aligned}
\end{equation}
where we have used (\ref{e421}) and (\ref{e422}). Integrating (\ref{e423}) w.r.t. time we have
\begin{equation*}
\begin{aligned}
\int_{t_1}^{t_2}\|\Delta P_m(\phi_m^3)-\Delta (\phi^3)\|^2ds&\le C\big(\|\phi_m-\phi\|_{L^2}+\|(P_m-I)\phi^3\|_{L^2}\big)\int_{t_1}^{t_2}\big(\|\Delta^2\phi_m\|_{L^2}+\|\Delta^2\phi\|_{L^2}\big)ds\\
&\le C\big(\|\phi_m-\phi\|_{L^2}+\|(P_m-I)\phi^3\|_{L^2}\big),
\end{aligned}
\end{equation*}
which, together with the strong convergence of $\phi_m$ in $C([0,T];L^2)$ and the convergence $\|(P_m-I)\phi^3\|_{L^2}\to 0$, implies that
\begin{equation}\label{e424}
\begin{aligned}
\lim_{m\to\infty}\int_{t_1}^{t_2}\|\Delta P_m(\phi_m^3)-\Delta (\phi^3)\|^2ds\to 0,\ \ \ \forall\ t_1,t_2\in[0,T].
\end{aligned}
\end{equation}
By H\"{o}lder's inequality we have
\begin{equation}\label{e425}
\begin{aligned}
\bigg|\int_{t_1}^{t_2}(\Delta P_m(\phi_m^3),\omega)ds-\int_{t_1}^{t_2}(\Delta (\phi^3),\omega)ds\bigg|
\le \Big(\int_{t_1}^{t_2}\|\Delta P_m(\phi_m^3)-\Delta (\phi^3)\|^2ds\Big)^{\frac{1}{2}}\|\omega\|T^{\frac{1}{2}}.
\end{aligned}
\end{equation}
Hence, by passing to the limit in $m$ and using (\ref{e424}) we have
\begin{equation}\label{e426}
\lim_{m\to\infty}\bigg|\int_{t_1}^{t_2}(\Delta P_m(\phi_m^3),\omega)ds-\int_{t_1}^{t_2}(\Delta (\phi^3),\omega)ds\bigg|\to 0,\ \ \ \forall\ t_1,t_2\in[0,T].
\end{equation}
This gives (\ref{e420}). As a consequence of (\ref{e419''}) and (\ref{e420}) we then have
\begin{equation}\label{e427}
\int_{t_1}^{t_2}(\Delta \mu_m,\omega)ds-\int_{t_1}^{t_2}(\Delta \mu,\omega)ds \longrightarrow 0,
\end{equation}
as $m\to\infty$, where $\mu:=\phi^3-\phi-\Delta\phi$.

Upon taking the $L^2$ inner product of $(\ref{ns})_1$ with $\omega\in H$ and integrating over the interval $[t_1,t_2]\subset [0,T]$ we obtain
\begin{equation}\label{e428}
(\phi_m(t_2),\omega)-(\phi_m(t_1),\omega)+\int_{t_2}^{t_1}(P_m(\vec{\mathbf{v}}_m\cdot\nabla\phi_m),\omega)ds
=\int_{t_2}^{t_1}(\Delta\mu_m,\omega)ds.
\end{equation}

In view of (\ref{e419'}) and (\ref{e427}) we see that it remains to show
\begin{equation}\label{e429}
\int_{t_1}^{t_2}(P_m(\vec{\mathbf{v}}_m\cdot\nabla\phi_m),\omega) ds-\int_{t_1}^{t_2}((\vec{\mathbf{v}}\cdot\nabla\phi),\omega) ds \longrightarrow 0,
\end{equation}
as $m\to\infty$ for every $t_1,t_2\in[0,T]$ and every $\omega\in H$, in order to prove (\ref{ss}). First, we notice that for every $\omega\in H$,
\begin{equation}\label{e430}
\begin{aligned}
\int_0^T\|(\nabla\phi)\omega\|^2ds & \le \int_0^T\|\nabla\phi\|_{L^\infty}^2\|\omega\|^2ds\\
&\le C\|\omega\|^2\int_0^T\|\phi\|_{H^3}^2ds\le C.
\end{aligned}
\end{equation}
This implies that $(\nabla\phi)\omega\in L^2([0,T];(L^2(\Omega))^2)$. Moreover, by H\"{o}lder's inequality we have
\begin{equation}\label{e431}
\begin{aligned}
&\bigg|\int_{t_1}^{t_2}(P_m(\vec{\mathbf{v}}_m\cdot\nabla\phi_m),\omega) ds-\int_{t_1}^{t_2}((\vec{\mathbf{v}}\cdot\nabla\phi),\omega) ds\bigg|\\
=&\bigg|\int_{t_1}^{t_2}(P_m(\vec{\mathbf{v}}_m\cdot\nabla\phi_m)-P_m(\vec{\mathbf{v}}_m\cdot\nabla\phi),\omega) ds+\int_{t_1}^{t_2}(P_m(\vec{\mathbf{v}}_m\cdot\nabla\phi)-P_m(\vec{\mathbf{v}}\cdot\nabla\phi),\omega) ds\bigg|\\
&+\int_{t_1}^{t_2}(P_m(\vec{\mathbf{v}}\cdot\nabla\phi)-(\vec{\mathbf{v}}\cdot\nabla\phi),\omega) ds\bigg|\\
\le &\int_{t_1}^{t_2}\|\vec{\mathbf{v}}_m\cdot\nabla(\phi_m-\phi)\|\|\omega\| ds+\int_{t_1}^{t_2}\|(\vec{\mathbf{v}}_m-\vec{\mathbf{v}})\cdot\nabla\phi\|\|\omega\| ds\\
&+\int_{t_1}^{t_2}\|P_m(\vec{\mathbf{v}}\cdot\nabla\phi)-(\vec{\mathbf{v}}\cdot\nabla\phi)\|\|\omega\| ds\\
\le &\|\omega\|\int_{t_1}^{t_2}\|\vec{\mathbf{v}}_m\|_{H^1}\|\nabla(\phi_m-\phi)\|_{H^1} ds+\|\omega\|\int_{t_1}^{t_2}\|\vec{\mathbf{v}}_m-\vec{\mathbf{v}}\|\|\nabla\phi\|_{L^\infty} ds\\
&+\|\omega\|\int_{t_1}^{t_2}\|P_m(\vec{\mathbf{v}}\cdot\nabla\phi)-(\vec{\mathbf{v}}\cdot\nabla\phi)\| ds\\
\le & \|\omega\|\bigg(\int_{t_1}^{t_2}\|\vec{\mathbf{v}}_m\|_{H^1}^2ds\bigg)^{\frac{1}{2}}\bigg(\int_{t_1}^{t_2}\|\phi_m-\phi\|_{H^2}^2 ds\bigg)^{\frac{1}{2}}\\
&+\|\omega\|\bigg(\int_{t_1}^{t_2}\|\vec{\mathbf{v}}_m-\vec{\mathbf{v}}\|^2ds\bigg)^{\frac{1}{2}}\bigg(\int_{t_1}^{t_2}\|\phi\|_{H^3}^2 ds\bigg)^{\frac{1}{2}}\\
&+\|\omega\|\bigg(\int_{t_1}^{t_2}\|P_m(\vec{\mathbf{v}}\cdot\nabla\phi)-(\vec{\mathbf{v}}\cdot\nabla\phi)\|^2ds\bigg)
^{\frac{1}{2}}|t_1-t_2|^{\frac{1}{2}}.
\end{aligned}
\end{equation}
Since $\vec{\mathbf{v}}_m$ is bounded in $L^2([0,T];\mathrm{\mathbf{V}})$ and the sequence $\{\phi_m\}$ converges strongly to $\phi$ in
$L^2([0,T];H^3)$, we know that the first term on the right hand side of (\ref{e431}) vanishes as $m\to\infty$. For the second term, since
$\vec{\mathbf{v}}_m$ converges strongly to $\vec{\mathbf{v}}$ in $L^2([0,T];\mathrm{\mathbf{H}})\subset L^2([0,T];(L^2(\Omega))^2)$, we know that it also goes to zero as $m\to\infty$. Furthermore, the last term tends to zero as $m\to\infty$ by virtue of the pointwise convergence of the projection operator. This proves (\ref{e429}). Thus by letting $m\to\infty$ in (\ref{e428}) and using (\ref{e419'}), (\ref{e427}) and (\ref{e429}) we get
\begin{equation}
(\phi(t_2),\omega)-(\phi(t_1),\omega)+\int_{t_2}^{t_1}((\vec{\mathbf{v}}\cdot\nabla\phi),\omega)ds
=\int_{t_2}^{t_1}(\Delta\mu,\omega)ds.
\end{equation}
Moreover, by using the strong convergence of $\{\phi_m\}$ in $L^2([0,T];H^3)$ and the estimate (\ref{e419}) one can show that $\vec{\mathbf{v}}_m=\mathrm{\mathbf{P}}_m\mathrm{\mathbf{P}}_\sigma(\mu_m\nabla\phi_m)$ also converges weakly to $\mathrm{\mathbf{P}}_\sigma(\mu\nabla\phi)$ in $L^2([0,T];\mathrm{\mathbf{H}})$. This gives that $\vec{\mathbf{v}}=\mathrm{\mathbf{P}}_\sigma(\mu\nabla\phi)$.

In order to complete the proof, it remains to show that $\phi\in C([0,T];H^2)$. For this purpose, we first re-visit the estimate (\ref{e412}), which reads
\begin{equation}
\begin{aligned}
\frac{1}{2}\frac{d}{dt}\|\Delta\phi_m\|^2+\frac{1}{2}\|\Delta^2\phi_m\|^2
\le C\|\nabla\Delta\phi_m\|^2\|\Delta\phi_m\|^2+C\|\nabla\Delta\phi_m\|^2.
\end{aligned}\label{e432}
\end{equation}
We observe that, by virtue of the boundary condition $\nabla\phi_m\cdot\mathbf{n}|_{\partial\Omega}=0$
$$
\|\nabla\Delta\phi_m\|^2=-(\Delta\phi_m,\Delta^2\phi_m)\le \|\Delta\phi_m\|\|\Delta^2\phi_m\|,
$$
which implies, by Young's inequality, that
\begin{equation}
\begin{aligned}
\frac{1}{2}\frac{d}{dt}\|\Delta\phi_m\|^2+\frac{1}{2}\|\Delta^2\phi_m\|^2 & \le C\|\Delta^2\phi_m\|\|\Delta\phi_m\|^3+C\|\Delta\phi_m\|\|\Delta^2\phi_m\|\\
&\le \frac{1}{4}\|\Delta^2\phi_m\|^2+C\big(\|\Delta\phi_m\|^6+\|\Delta\phi_m\|^2\big).
\end{aligned}
\end{equation}
After rearranging terms we have
\begin{equation}
\begin{aligned}
2\frac{d}{dt}\|\Delta\phi_m\|^2+\|\Delta^2\phi_m\|^2 \le C\big(\|\Delta\phi_m\|^6+\|\Delta\phi_m\|^2\big)\le C,
\end{aligned}\label{e432'}
\end{equation}
where, in the last step, we used (\ref{e419}). Integrating (\ref{e432'}) over $[s,t]\subset [0,T]$ for $s<t$ we have
\begin{equation}\label{e433}
\begin{aligned}
\|\Delta\phi_m(t)\|^2+\int_s^t\|\Delta^2\phi_m(\tau)\|^2d\tau
\le \|\Delta\phi_m(s)\|^2+C(t-s).
\end{aligned}
\end{equation}
Due to (\ref{E}) and the inequality,
$$
\liminf_{n\to\infty}(a_n)+\liminf_{n\to\infty}(b_n)\le \limsup_{n\to\infty}(a_n+b_n)
$$
for any sequences $a_n,b_n\ge0$, we have
\begin{equation}\label{e434}
\begin{aligned}
\|\Delta\phi(t)\|^2+\int_s^t\|\Delta^2\phi(\tau)\|^2d\tau
\le \|\Delta\phi(s)\|^2+C(t-s),
\end{aligned}
\end{equation}
for every $s,t\in E$, $s<t$, $E$ is given in (\ref{E}). After rearranging terms in (\ref{e434}) we have
\begin{equation}\label{e435}
\begin{aligned}
\big|\|\Delta\phi(t)\|^2-\|\Delta\phi(s)\|^2\big|\le \int_s^t\|\Delta^2\phi(\tau)\|^2d\tau+C(t-s),
\end{aligned}
\end{equation}
which can be extended to all $s$ and $t$ in $[0,T]$ by manipulating sequences in $E$ that converge to elements in $[0,T]\backslash E$. Using (\ref{e435}) and the fact that $\phi\in C_W([0,T];H^2)$ we get $\phi\in C([0,T];H^2)$. Finally, we notice that $T>0$ is arbitrary. This gives the globally existence in 2D. This completes the proof of Theorem 1.1.

\section{Local existence in 3D}

In this section we prove Theorem 1.2. Now we turn to (\ref{e330}) in 3D. Instead of applying the Cauchy-Schwarz inequality to all the terms on the RHS of (\ref{e330}), we keep the
first term and apply the Cauchy-Schwarz and Young's inequalities to the other two terms to get
\begin{equation}
\begin{aligned}
&\frac{1}{2}\frac{d}{dt}\|\Delta\phi_m\|^2+\|\Delta^2\phi_m\|^2\\=&-\int_\Omega (\vec{\mathbf{v}}_m\cdot\nabla\phi_m)(\Delta^2\phi_m) d\mathbf{x}+
\int_\Omega\Delta((\phi_m+\bar{\phi})^3)(\Delta^2\phi_m)d\mathbf{x}+\int_\Omega\Delta(\phi_m)(\Delta^2\phi_m)d\mathbf{x}\\
\le &\frac{1}{4}\|\Delta^2\phi_m\|^2-\int_\Omega (\vec{\mathbf{v}}_m\cdot\nabla\phi_m)(\Delta^2\phi_m) d\mathbf{x}+2\|\Delta((\phi_m+\bar{\phi})^3)\|^2+2\|\Delta\phi_m\|^2.
\end{aligned}\label{e51}
\end{equation}

By applying Lemma 2.2 (IX), we estimate the second term on the RHS of (\ref{e51}) as:
\begin{equation}
\begin{aligned}
&\bigg|-\int_\Omega (\vec{\mathbf{v}}_m\cdot\nabla\phi_m)(\Delta^2\phi_m) d\mathbf{x}\bigg|\\
\le &\|(\vec{\mathbf{v}}_m\cdot\nabla\phi_m)\|\|\Delta^2\phi_m\|\\
\le &\|\vec{\mathbf{v}}_m\|_{L^3}\|\nabla\phi_m\|_{L^6}\|\Delta^2\phi_m\|\\
\le &C\big(\|\vec{\mathbf{v}}_m\|^{1/2}\|\nabla \vec{\mathbf{v}}_m\|^{1/2}+\|\vec{\mathbf{v}}_m\|\big)\|\nabla\phi_m\|_{H^1}\|\Delta^2\phi_m\|\\
\le &C\|\mu_m\nabla\phi_m\|^{1/2}\big(\|\nabla\mu_m\cdot\nabla\phi_m\|^{1/2}+\|\mu_m\nabla^2\phi_m\|^{1/2}+\|\mu_m\nabla\phi_m\|^{1/2}\big)
\|\Delta\phi_m\|\|\Delta^2\phi_m\|,
\end{aligned}\label{e52}
\end{equation}
where in the last inequality we used $(\ref{ns})_3$ and Lemma 2.2 (XII).

For the first term on the RHS of (\ref{e52}), by Lemma 2.2 (X) and (XII) and (\ref{e43}) we have
\begin{equation}
\begin{aligned}
\|\mu_m\nabla\phi_m\|^{1/2}\le &\|\mu_m\|^{1/2}_{L^3}\|\nabla\phi_m\|^{1/2}_{L^6}\\
\le &C\|\mu_m\|^{1/4}\|\nabla\mu_m\|^{1/4}\|\nabla\phi_m\|_{H^1}^{1/2}\\
\le &C\|\Delta\phi_m\|^{1/4}\|\Delta\mu_m\|^{1/4}\|\Delta\phi_m\|^{1/2}\\
\le &C\|\Delta\mu_m\|^{1/4}\|\Delta\phi_m\|^{3/4}.
\end{aligned}\label{e53}
\end{equation}
In order to estimate $\|\Delta\mu_m\|^{1/4}$, by (\ref{e414}) we have
\begin{equation}
\|\Delta\mu_m\|^{1/4}\le C\big(\|\Delta\phi_m\|^{1/2}+\|\Delta^2\phi_m\|^{1/4}\big).\label{e54}
\end{equation}
Plugging (\ref{e54}) into (\ref{e53}) we have
\begin{equation}
\|\mu_m\nabla\phi_m\|^{1/2}\le C\big(\|\Delta\phi_m\|^{1/2}+\|\Delta^2\phi_m\|^{1/4}\big)\|\Delta\phi_m\|^{3/4}.\label{e55}
\end{equation}

Next, to estimate the three terms in the bracket on the RHS of (\ref{e52}), we use similar arguments and Lemma 2.2 (XII) and (XIII) to obtain
\begin{equation}
\begin{aligned}
\|\nabla\mu_m\cdot\nabla\phi_m\|^{1/2}&\le \|\nabla\mu_m\|_{L^3}^{1/2}\|\nabla\phi_m\|_{L^6}^{1/2}\\
&\le \|\nabla\mu_m\|_{H^1}^{1/2}\|\nabla\phi_m\|_{H^1}^{1/2}\\
&\le C\|\Delta\mu_m\|^{1/2}\|\Delta\phi_m\|^{1/2}\\
&\le C\big(\|\Delta\phi_m\|+\|\Delta^2\phi_m\|^{1/2}\big)\|\Delta\phi_m\|^{1/2},
\end{aligned}\label{e56}
\end{equation}
\begin{equation}
\begin{aligned}
\|\mu_m\nabla^2\phi_m\|^{1/2}&\le \|\mu_m\|_{L^\infty}^{1/2}\|\Delta\phi_m\|^{1/2}\\
&\le C\|\mu_m\|_{H^2}^{1/2}\|\Delta\phi_m\|^{1/2}\\
&\le C\|\Delta\mu_m\|^{1/2}\|\Delta\phi_m\|^{1/2}\\
&\le C\big(\|\Delta\phi_m\|+\|\Delta^2\phi_m\|^{1/2}\big)\|\Delta\phi_m\|^{1/2},
\end{aligned}\label{e57}
\end{equation}
and
\begin{equation}
\begin{aligned}
\|\mu_m\nabla\phi_m\|^{1/2}&\le \|\mu_m\|_{L^3}^{1/2}\|\nabla\phi_m\|_{L^6}^{1/2}\\
&\le \|\mu_m\|_{H^1}^{1/2}\|\nabla\phi_m\|_{H^1}^{1/2}\\
&\le C\|\Delta\mu_m\|^{1/2}\|\Delta\phi_m\|^{1/2}\\
&\le C\big(\|\Delta\phi_m\|+\|\Delta^2\phi_m\|^{1/2}\big)\|\Delta\phi_m\|^{1/2}.
\end{aligned}\label{e58}
\end{equation}
Therefore, we have
\begin{equation}
\|\nabla\mu_m\cdot\nabla\phi_m\|^{1/2}+\|\mu_m\nabla^2\phi_m\|^{1/2}+\|\mu_m\nabla\phi_m\|^{1/2}\le \big(\|\Delta\phi_m\|+\|\Delta^2\phi_m\|^{1/2}\big)\|\Delta\phi_m\|^{1/2}.\label{e59}
\end{equation}

Combining (\ref{e55}) and (\ref{e59}) we get
\begin{equation}
\begin{aligned}
&\|\mu_m\nabla\phi_m\|^{1/2}\big(\|\nabla\mu_m\cdot\nabla\phi_m\|^{1/2}+\|\mu_m\nabla^2\phi_m\|^{1/2}+\|\mu_m\nabla\phi_m\|^{1/2}\big)\\
\le &C\big(\|\Delta\phi_m\|^{3/2}+\|\Delta^2\phi_m\|^{3/4}\big)\|\Delta\phi_m\|^{5/4},\label{e510}
\end{aligned}
\end{equation}
which, together with (\ref{e52}), implies that
\begin{equation}
\bigg|-\int_\Omega (\vec{\mathbf{v}}_m\cdot\nabla\phi_m)(\Delta^2\phi_m) d\mathbf{x}\bigg|\le C\|\Delta\phi_m\|^{9/4}\|\Delta^2\phi_m\|^{7/4}+C\|\Delta\phi_m\|^{15/4}\|\Delta^2\phi_m\|.\label{e511}
\end{equation}

The estimate for the second term on the RHS of (\ref{e51}) follows from (\ref{e335}):
\begin{equation}
\|\Delta((\phi_m+\bar{\phi})^3)\|^2\le C\big(\|\Delta\phi_m\|^4+\|\Delta\phi_m\|^2\big).\label{e512}
\end{equation}

Plugging the estimates (\ref{e511}) and (\ref{e512}) into (\ref{e51}) and applying Young's inequality we get
\begin{equation}
\begin{aligned}
&\frac{1}{2}\frac{d}{dt}\|\Delta\phi_m\|^2+\|\Delta^2\phi_m\|^2\\
\le &\frac{1}{4}\|\Delta^2\phi_m\|^2+
C\|\Delta\phi_m\|^{9/4}\|\Delta^2\phi_m\|^{7/4}+\\
&C\|\Delta\phi_m\|^{15/4}\|\Delta^2\phi_m\|+C\big(\|\Delta\phi_m\|^4+\|\Delta\phi_m\|^2\big)\\
\le &\frac{1}{2}\|\Delta^2\phi_m\|^2+C\big(\|\Delta\phi_m\|^{18}+\|\Delta\phi_m\|^2\big),
\end{aligned}\label{e513}
\end{equation}
which implies that
\begin{equation}
\frac{d}{dt}\|\Delta\phi_m\|^2+\|\Delta^2\phi_m\|^2
\le C\big(\|\Delta\phi_m\|^{18}+\|\Delta\phi_m\|^2\big).\label{e514}
\end{equation}
We remark that the constant $C$ on the RHS of (\ref{e514}) is independent of time.

Let
\begin{equation}
y_m(t)=\|\Delta\phi_m(t)\|^2+1.\label{e515}
\end{equation}
Then from (\ref{e514}) we have
\begin{equation}
y_m'(t)\le C(y_m(t))^9,\label{e516}
\end{equation}
which implies immediately that
\begin{equation}
y_m(t)\le \frac{y_m(0)}{(1-8y_m(0)^8Ct)^{1/8}},\ \ \forall\ 0\le t<\frac{1}{8y_m(0)^8C}\equiv T_m.\label{e517}
\end{equation}
Thus for any $m$ we have
\begin{equation}
y_m(t)\le \frac{\|\Delta\phi^0\|^2+1}{(1/2)^{1/8}},\ \ \forall\ 0\le t\le \frac{1}{16(\|\Delta\phi^0\|^2+1)^8C}\equiv T^*,\label{e518}
\end{equation}
since $y_m(0)\le y(0)=\|\Delta\phi^0\|^2+1$.

As consequences of above estimates we have
\begin{equation}
\sup_{\tau\in[0,T^*]}\big(\|\phi_m(\tau)\|^2_{H^2}+\|\mu_m(\tau)\|^2\big)+\int_0^{T^*}\big(\|\vec{\mathbf{v}}_m\|^2_{H^1}+
\|\phi_m\|^2_{H^4}+\|\mu_m\|^2_{H^2}+\Big\|\frac{d\phi_m}{dt}\Big\|^2\big)ds\le C.\label{e519}
\end{equation}

These estimates, combined with the compactness arguments given in the preceding section, imply the local existence result in 3D. This completes the
proof of Theorem 1.2.

\section{Uniqueness and continuous dependence in initial data}

In this section, we prove the uniqueness and continuous dependence on initial data of the solutions in both 2D and 3D. Suppose that there are two solutions, namely, $(\phi_1,\mu_1,\vec{\mathbf{v}}_1,P_1)$ and $(\phi_2,\mu_2,\vec{\mathbf{v}}_2,P_2)$ on a joint interval of existence $[0,\bar{T}]$. We use $\tilde{f}$ to denote the difference of two functions. Then the differences satisfy
\begin{equation}
\left\{
\begin{aligned}
&\tilde{\phi}_t+\tilde{\vec{\mathbf{v}}}\cdot\nabla\phi_1+\vec{\mathbf{v}}_2\cdot\nabla\tilde{\phi}=\Delta\tilde{\mu},\\
&\tilde{\mu}=\phi_1^3-\phi_2^3-\tilde{\phi}-\Delta\tilde{\phi},\\
&\tilde{\vec{\mathbf{v}}}=-\nabla\tilde{P}+\tilde{\mu}\nabla\phi_1+\mu_2\nabla\tilde{\phi};\\
&\nabla\tilde{\phi}\cdot\mathbf{n}|_{\partial\Omega}=\nabla\tilde{\mu}\cdot\mathbf{n}|_{\partial\Omega}=\tilde{\vec{\mathbf{v}}}\cdot\mathbf{n}|_{\partial\Omega}=0,\\
&\tilde{\phi}(\mathbf{x},0)=\phi_1^0(\mathbf{x})-\phi_2^0(\mathbf{x}).
\end{aligned}
\right.\label{e61}
\end{equation}

Taking the $L^2$ inner products of $(\ref{e61})_1$ with $\tilde{\mu}$ and $(\ref{e61})_3$ with $\tilde{\vec{\mathbf{v}}}$, respectively, and adding the results we have
\begin{equation}
\begin{aligned}
\frac{d}{dt}\bigg(\frac{1}{2}\|\nabla\tilde{\phi}\|^2-\frac{1}{2}\|\tilde{\phi}\|^2\bigg)+\int_\Omega(\phi^3_1-\phi_2^3)\frac{d\tilde{\phi}}{dt}d\mathbf{x}
+\|\tilde{\vec{\mathbf{v}}}\|^2+\|\nabla\tilde{\mu}\|^2&=\int_\Omega\mu_2\nabla\tilde{\phi}\cdot\tilde{\vec{\mathbf{v}}}d\mathbf{x}-\int_\Omega \tilde{\mu}\vec{\mathbf{v}}_2\cdot\nabla\tilde{\phi}d\mathbf{x}\\
&=\int_\Omega\mu_2\nabla\tilde{\phi}\cdot\tilde{\vec{\mathbf{v}}}d\mathbf{x}+\int_\Omega \tilde{\phi}\vec{\mathbf{v}}_2\cdot\nabla\tilde{\mu}d\mathbf{x}\\
\end{aligned}\label{e62}
\end{equation}

Define
\begin{equation}
G(\tilde{\phi})=\int_0^{\tilde{\phi}}\big((\phi_2+\xi)^3-\phi_2^3\big)d\xi.\label{e63}
\end{equation}
Then it is easy to see that
\begin{equation}
G(\tilde{\phi})\ge0\ \ \ \mathrm{for}\ \tilde{\phi}\in\mathbb{R},\ \ G(\tilde{\phi})=0\ \ \mathrm{iff}\ \ \tilde{\phi}=0,\ \ \mathrm{and}\ \ G'(\tilde{\phi})=(\phi^3_1-\phi_2^3).\label{e64}
\end{equation}
Since, by the existence theorems, $\frac{d\tilde{\phi}}{dt}\in L^2([0,\bar{T}];H)$, then $\tilde{\phi}\in C([0,\bar{T}];H)$ and is an absolutely continuous function. Moreover, since $G(\tilde{\phi})$ is a polynomial, it is differentiable, and $G(\tilde{\phi})$ is absolutely continuous in time with value in $H$ and $\frac{dG(\tilde{\phi})}{dt}=G'(\tilde{\phi})\frac{d\tilde{\phi}}{dt}\in L^2([0,\bar{T}];H)$. Furthermore, one can also show that (c.f. \cite{T})
\begin{equation}
\frac{d}{dt}\int_\Omega G(\tilde{\phi})d\mathbf{x}=\int_\Omega G'(\tilde{\phi})\frac{d\tilde{\phi}}{dt}d\mathbf{x}.\label{e65}
\end{equation}

Then plugging (\ref{e65}) into (\ref{e62}) and using (\ref{e64}) we have
\begin{equation}
\frac{d}{dt}\bigg(\frac{1}{2}\|\nabla\tilde{\phi}\|^2-\frac{1}{2}\|\tilde{\phi}\|^2+\int_\Omega G(\tilde{\phi})d\mathbf{x}\bigg)
+\|\tilde{\vec{\mathbf{v}}}\|^2+\|\nabla\tilde{\mu}\|^2\le \|\mu_2\nabla\tilde{\phi}\|^2+\frac{1}{4}\|\tilde{\vec{\mathbf{v}}}\|^2+\frac{1}{4}\|\nabla\tilde{\mu}\|^2+\|\tilde{\phi}\vec{\mathbf{v}}_2\|^2.\label{e66}
\end{equation}

We estimate the two nonlinear terms on the RHS of (\ref{e66}) as:
\begin{equation}
\begin{aligned}
\|\mu_2\nabla\tilde{\phi}\|^2+\|\tilde{\phi}\vec{\mathbf{v}}_2\|^2&\le \|\mu_2\|_{L^\infty}^2\|\nabla\tilde{\phi}\|^2+\|\tilde{\phi}\|_{L^6}^2\|\vec{\mathbf{v}}_2\|_{L^3}^2\\
&\le C\|\mu_2\|_{H^2}^2\|\nabla\tilde{\phi}\|^2+C\big(\|\nabla\tilde{\phi}\|^2+\|\tilde{\phi}\|^2\big)\|\vec{\mathbf{v}}_2\|_{H^1}^2,
\end{aligned}\label{e67}
\end{equation}
where we have used Poincar\'{e}'s inequality since $\phi_1$ and $\phi_2$ have the same average over $\Omega$. Then we get
\begin{equation}
\frac{d}{dt}\bigg(\frac{1}{2}\|\nabla\tilde{\phi}\|^2-\frac{1}{2}\|\tilde{\phi}\|^2+\int_\Omega G(\tilde{\phi})d\mathbf{x}\bigg)
+\frac{3}{4}\|\tilde{\vec{\mathbf{v}}}\|^2+\frac{3}{4}\|\nabla\tilde{\mu}\|^2\le C\big(\|\mu_2\|_{H^2}^2+\|\vec{\mathbf{v}}_2\|_{H^1}^2\big)\big(\|\nabla\tilde{\phi}\|^2+\|\tilde{\phi}\|^2\big).\label{e68}
\end{equation}

Next, taking the $L^2$ inner product of $(\ref{e61})_1$ with $\tilde{\phi}$ we have
\begin{equation}
\begin{aligned}
\frac{1}{2}\frac{d}{dt}\|\tilde{\phi}\|^2&=-\int_\Omega (\tilde{U}\cdot\nabla\phi_1)\tilde{\phi}d\mathbf{x}+\int_\Omega\tilde{\mu}\Delta\tilde{\phi}d\mathbf{x}\\
&=\int_\Omega (\tilde{\vec{\mathbf{v}}}\cdot\nabla\tilde{\phi}){\phi_1}d\mathbf{x}+\int_\Omega\tilde{\mu}\Delta\tilde{\phi}d\mathbf{x}\\
&=\int_\Omega (\tilde{\vec{\mathbf{v}}}\cdot\nabla\tilde{\phi}){\phi_1}d\mathbf{x}-\|\Delta\tilde{\phi}\|^2+\|\nabla\tilde{\phi}\|^2+\int_\Omega(\phi_1^2+\phi_1\phi_2+\phi_2^2)\tilde{\phi}\Delta\tilde{\phi}d\mathbf{x}\\
&\le \|\tilde{\vec{\mathbf{v}}}\|\|\phi_1\|_{L^\infty}\|\nabla\tilde{\phi}\|-\|\Delta\tilde{\phi}\|^2+\|\nabla\tilde{\phi}\|^2+\|\phi_1^2+\phi_1\phi_2+\phi_2^2\|_{L^\infty}\|\tilde{\phi}\|\|\Delta\tilde{\phi}\|\\
&\le \frac{1}{8}\|\tilde{\vec{\mathbf{v}}}\|^2-\frac{1}{4}\|\Delta\tilde{\phi}\|^2+C\big(\|\nabla\tilde{\phi}\|^2+\|\tilde{\phi}\|^2\big),
\end{aligned}\label{e69}
\end{equation}
which gives
\begin{equation}
\frac{d}{dt}\|\tilde{\phi}\|^2+\frac{1}{2}\|\Delta\tilde{\phi}\|^2\le \frac{1}{4}\|\tilde{\vec{\mathbf{v}}}\|^2+C\big(\|\nabla\tilde{\phi}\|^2+\|\tilde{\phi}\|^2\big).\label{e610}
\end{equation}

Combining (\ref{e68}) and (\ref{e610}) and since $G(\tilde{\phi})\ge0$ we have
\begin{equation}
\frac{d}{dt}\bigg(\frac{1}{2}\|\nabla\tilde{\phi}\|^2+\frac{1}{2}\|\tilde{\phi}\|^2+\int_\Omega G(\tilde{\phi})d\mathbf{x}\bigg)
+\frac{1}{2}\|\tilde{\vec{\mathbf{v}}}\|^2+\frac{3}{4}\|\nabla\tilde{\mu}\|^2\le C(t)\bigg(\|\nabla\tilde{\phi}\|^2+\|\tilde{\phi}\|^2+\int_\Omega G(\tilde{\phi})d\mathbf{x}\bigg),\label{e611}
\end{equation}
where
\begin{equation}
C(t)=C\big(\|\mu_2\|_{H^2}^2+\|\vec{\mathbf{v}}_2\|_{H^1}^2+1\big)\ \ \ \mathrm{satisfying}\ \int_0^{\bar{T}} C(\tau)ds\le C,\ \ \forall\ t\in[0,\bar{T}].\label{e612}
\end{equation}
Applying Gronwall's inequality to (\ref{e611}) we get immediately that for any $t\in[0,\bar{T}]$,
\begin{equation}\label{e613}
\begin{aligned}
&\bigg(\frac{1}{2}\|\nabla\tilde{\phi}(t)\|^2+\frac{1}{2}\|\tilde{\phi}(t)\|^2+\int_\Omega G(\tilde{\phi}(t))d\mathbf{x}\bigg)\\
\le &\exp\bigg\{\int_0^{\bar{T}}C(\tau)ds\bigg\}
\bigg(\frac{1}{2}\|\nabla\tilde{\phi}^0\|^2+\frac{1}{2}\|\tilde{\phi}^0\|^2+\int_\Omega G(\tilde{\phi}^0)d\mathbf{x}\bigg).
\end{aligned}
\end{equation}
Therefore, the solutions depend continuously on the initial data in the sense of (\ref{e613}). In particular, if $\tilde{\phi}^0=0$ we obtain the uniqueness of the solutions. This completes the proof of Theorem 1.3.

\section{Higher order spatial regularity}

In this section, we show that the unique solution obtained in previous sections enjoys higher order spatial regularity.
The idea is adopted from the proof for the Navier-Stokes equations in \cite{CF,T}.

Taking the $L^2$ inner product of $(\ref{e31})_1$ with $\Delta^3\phi_m$ we have
\begin{equation}
\begin{aligned}
&\frac{1}{2}\frac{d}{dt}\|\nabla\Delta\phi_m\|^2+\|\nabla\Delta^2\phi_m\|^2\\
=&\int_\Omega \Delta^3\phi_m\big(P_m(\vec{\mathbf{v}}_m\cdot\nabla\phi_m)-\Delta(\phi_m+\bar{\phi})^3+\Delta\phi_m\big)d\mathbf{x}\\
=&-\int_\Omega \nabla\Delta^2\phi_m\cdot\nabla\big(P_m(\vec{\mathbf{v}}_m\cdot\nabla\phi_m)-\Delta(\phi_m+\bar{\phi})^3+\Delta\phi_m\big)d\mathbf{x}\\
\le &\frac{1}{2}\|\nabla\Delta^2\phi_m\|^2+\frac{3}{2}\big(\|\nabla\Delta\phi_m\|^2+\|\nabla\Delta((\phi_m+\bar{\phi})^3)\|^2+
\|\nabla(\vec{\mathbf{v}}_m\cdot\nabla\phi_m)\|^2\big).
\end{aligned}\label{e71}
\end{equation}
We break the estimate of the RHS of (\ref{e71}) into several steps.

{\bf Step 1.} For the third term on the RHS of (\ref{e71}), after direct calculations we have
\begin{equation}
\begin{aligned}
&\|\nabla\Delta((\phi_m+\bar{\phi})^3)\|^2\\
\le &C\big(\|\nabla\phi_m\|^6_{L^6}+\|(\phi_m+\bar{\phi})\|_{L^\infty}^2\|\nabla\phi_m\|_{L^4}^2\|\nabla^2\phi_m\|_{L^4}^2
+\|(\phi_m+\bar{\phi})\|_{L^\infty}^4\|\nabla\Delta\phi_m\|^2\big).
\end{aligned}\label{e72}
\end{equation}
Using Sobolev embedding, Lemma 2.2 (XII) and (\ref{e321}) we estimate the first term on the RHS of (\ref{e72}) as
\begin{equation}
\|\nabla\phi_m\|^6_{L^6}\le C\|\nabla\phi_m\|^6_{H^1}\le C\|\Delta\phi_m\|^6\le C\|\Delta\phi_m\|^4\|\nabla\Delta\phi_m\|^2.\label{e73}
\end{equation}
Similarly, for the second and third terms, using the triangle inequality, Sobolev embeddings, Lemma 2.2 (XII) and (XIII) we have
\begin{equation}
\begin{aligned}
&C\big(\|(\phi_m+\bar{\phi})\|_{L^\infty}^2\|\nabla\phi_m\|_{L^4}^2\|\nabla^2\phi_m\|_{L^4}^2
+\|(\phi_m+\bar{\phi})\|_{L^\infty}^4\|\nabla\Delta\phi_m\|^2\big)\\
\le &C\big(\|\phi_m\|_{H^2}^2\|\nabla\phi_m\|_{H^1}^2\|\nabla^2\phi_m\|_{H^1}^2+\|\nabla\phi_m\|_{H^1}^2\|\nabla^2\phi_m\|_{H^1}^2
+\|\phi_m\|_{H^2}^4\|\nabla\Delta\phi_m\|^2+\|\nabla\Delta\phi_m\|^2\big)\\
\le &C\big(\|\Delta\phi_m\|^2\|\Delta\phi_m\|^2\|\nabla\Delta\phi_m\|^2+\|\Delta\phi_m\|^2\|\nabla\Delta\phi_m\|^2
+\|\Delta\phi_m\|^4\|\nabla\Delta\phi_m\|^2+\|\nabla\Delta\phi_m\|^2\big)\\
\le &C\big(\|\Delta\phi_m\|^4+1\big)\|\nabla\Delta\phi_m\|^2.
\end{aligned}\label{e74}
\end{equation}
Plugging (\ref{e73}) and (\ref{e74}) into (\ref{e72}) we have
\begin{equation}
\|\nabla\Delta((\phi_m+\bar{\phi})^3)\|^2\le C\big(\|\Delta\phi_m\|^4+1\big)\|\nabla\Delta\phi_m\|^2.\label{e75}
\end{equation}
Using previous estimates (\ref{e419}) and (\ref{e519}) we obtain for both 2D and 3D, within the interval of existence,
\begin{equation}
\|\nabla\Delta(\phi_m^3)\|^2\le C\|\nabla\Delta\phi_m\|^2.\label{e76}
\end{equation}

{\bf Step 2.} To estimate $\|\nabla(\vec{\mathbf{v}}_m\cdot\nabla\phi_m)\|^2$, first, we have
\begin{equation}
\|\nabla(\vec{\mathbf{v}}_m\cdot\nabla\phi_m)\|^2\le \|\vec{\mathbf{v}}_m\cdot\nabla^2\phi_m\|^2+\|\nabla(\vec{\mathbf{v}}_m)\cdot\nabla\phi_m\|^2.\label{e77}
\end{equation}
The first term on the RHS of (\ref{e77}) is estimated as
\begin{equation}
\begin{aligned}
\|\vec{\mathbf{v}}_m\cdot\nabla^2\phi_m\|^2&\le \|\vec{\mathbf{v}}_m\|_{L^3}^2\|\nabla^2\phi_m\|^2_{L^6}\\
&\le C\|\vec{\mathbf{v}}_m\|_{H^1}^2\|\nabla^2\phi_m\|^2_{H^1}\\
&\le C\|\vec{\mathbf{v}}_m\|_{H^1}^2\|\nabla\Delta\phi_m\|^2.
\end{aligned}\label{e78}
\end{equation}
For the second term, using Lemma 2.2 (XII) and (XIII) and (\ref{e321}) and (\ref{e322}) we have
\begin{equation}
\begin{aligned}
\|\nabla(\vec{\mathbf{v}}_m)\cdot\nabla\phi_m\|^2&\le \|\nabla(\vec{\mathbf{v}}_m)\|^2\|\nabla\phi_m\|^2_{L^\infty}\\
&\le C\|\nabla(\mu_m\nabla\phi_m)\|^2\|\nabla\phi_m\|^2_{H^2}\\
&\le C\big(\|\nabla\mu_m\cdot\nabla\phi_m\|^2+\|\mu_m\nabla^2\phi_m\|^2\big)\|\nabla\Delta\phi_m\|^2\\
&\le C\big(\|\nabla\mu_m\|_{L^4}^2\|\nabla\phi_m\|_{L^4}^2+\|\mu_m\|_{L^\infty}^2\|\nabla^2\phi_m\|^2\big)\|\nabla\Delta\phi_m\|^2\\
&\le C\big(\|\nabla\mu_m\|_{H^1}^2\|\nabla\phi_m\|_{H^1}^2+\|\mu_m\|_{H^2}^2\|\Delta\phi_m\|^2\big)\|\nabla\Delta\phi_m\|^2\\
&\le C\|\Delta\mu_m\|^2\|\Delta\phi_m\|^2\|\nabla\Delta\phi_m\|^2\\
&\le C\|\Delta\mu_m\|^2\|\nabla\Delta\phi_m\|^2,
\end{aligned}\label{e79}
\end{equation}
where we have used the uniform estimate of $\|\Delta\phi_m\|^2$.
Plugging (\ref{e78})--(\ref{e79}) into (\ref{e77}) we then have
\begin{equation}
\|\nabla(\vec{\mathbf{v}}_m\cdot\nabla\phi_m)\|^2\le C\big(\|\vec{\mathbf{v}}_m\|_{H^1}^2+\|\Delta\mu_m\|^2\big)\|\nabla\Delta\phi_m\|^2.\label{e710}
\end{equation}

Therefore, combining (\ref{e71}), (\ref{e76}) and (\ref{e710}) we have
\begin{equation}
\frac{d}{dt}\|\nabla\Delta\phi_m\|^2+\|\nabla\Delta^2\phi_m\|^2\le C\big(1+\|\vec{\mathbf{v}}_m\|_{H^1}^2+\|\Delta\mu_m\|^2\big)\|\nabla\Delta\phi_m\|^2.\label{e711}
\end{equation}

{\bf Step 3.} Taking the $L^2$ inner product of $(\ref{ns})_1$ with $\Delta^4\phi_m$ and using the Cauchy-Schwarz inequality we have
\begin{equation}
\begin{aligned}
\frac{1}{2}\frac{d}{dt}\|\Delta^2\phi_m\|^2+\frac{3}{4}\|\Delta^3\phi_m\|^2\le \Big(\|\Delta^2\phi_m\|^2+\|\Delta^2((\phi_m+\bar{\phi})^3)\|^2+
\|\Delta(\vec{\mathbf{v}}_m\cdot\nabla\phi_m)\|^2\Big).
\end{aligned}\label{e712}
\end{equation}
By a direct calculation, Sobolev embeddings and Lemma 2.2 (XII) and (XIII) we have
\begin{equation}
\begin{aligned}
&\|\Delta^2((\phi_m+\bar{\phi})^3)\|^2\\
\le C\big(&\|\nabla\phi_m\|^4_{L^4}\|\Delta\phi_m\|_{L^\infty}^2+\|(\phi_m+\bar{\phi})\|_{L^\infty}^2\|\Delta\phi_m\|_{L^\infty}^2\|\Delta\phi_m\|^2+\\
&\|(\phi_m+\bar{\phi})\|_{L^\infty}^2\|\nabla\phi_m\|_{L^3}^2\|\nabla\Delta\phi_m\|^2_{L^6}+\|(\phi_m+\bar{\phi})\|_{L^\infty}^4\|\Delta^2\phi_m\|^2\big)\\
\le C\big(&\|\Delta\phi_m\|^4\|\Delta^2\phi_m\|^2+\|\Delta\phi_m\|^2\|\Delta^2\phi_m\|^2\|\Delta\phi_m\|^2+\|\Delta^2\phi_m\|^2\|\Delta\phi_m\|^2+\\
&\|\Delta\phi_m\|^2\|\Delta\phi_m\|^2\|\Delta^2\phi_m\|^2+\|\Delta\phi_m\|^2\|\Delta^2\phi_m\|^2+\|\Delta\phi_m\|^4\|\Delta^2\phi_m\|^2+\|\Delta^2\phi_m\|^2\big)\\
\le C\big(&\|\Delta\phi_m\|^2+1\big)^2\|\Delta^2\phi_m\|^2,
\end{aligned}\label{e713}
\end{equation}
which implies, within the interval of existence, that
\begin{equation}
\begin{aligned}
\|\Delta^2((\phi_m+\bar{\phi})^3)\|^2\le C\|\Delta^2\phi_m\|^2,
\end{aligned}\label{e714}
\end{equation}
where we used (\ref{e419}) and (\ref{e519}) for the estimate of $\|\Delta\phi_m\|^2$.
So we update (\ref{e712}) as
\begin{equation}
\begin{aligned}
\frac{1}{2}\frac{d}{dt}\|\Delta^2\phi_m\|^2+\frac{3}{4}\|\Delta^3\phi_m\|^2\le C\|\Delta^2\phi_m\|^2+
\frac{3}{2}\|\Delta(\vec{\mathbf{v}}_m\cdot\nabla\phi_m)\|^2.
\end{aligned}\label{e715}
\end{equation}

{\bf Step 4.} Next, we estimate $\|\Delta(\vec{\mathbf{v}}_m\cdot\nabla\phi_m)\|^2$. First, we have
\begin{equation}
\|\Delta(\vec{\mathbf{v}}_m\cdot\nabla\phi_m)\|^2\le 4\big(\|\vec{\mathbf{v}}_m\cdot\nabla\Delta\phi_m\|^2+\|\nabla(\vec{\mathbf{v}}_m)\Delta\phi_m\|^2+\|\Delta(\vec{\mathbf{v}}_m)\nabla\phi_m\|^2\big).\label{e716}
\end{equation}

The first term on the RHS is estimated, using Lemma 2.2 (XII) and (XIII), as
\begin{equation}
\begin{aligned}
\|\vec{\mathbf{v}}_m\cdot\nabla\Delta\phi_m\|^2\le \|\vec{\mathbf{v}}_m\|_{L^3}^2\|\nabla\Delta\phi_m\|^2_{L^6}\le C\|\vec{\mathbf{v}}_m\|_{H^1}^2\|\nabla\Delta\phi_m\|^2_{H^1}\le C\|\vec{\mathbf{v}}_m\|_{H^1}^2\|\Delta^2\phi_m\|^2.
\end{aligned}\label{e717}
\end{equation}

Similarly, we have for the second term
\begin{equation}
\begin{aligned}
\|\nabla(\vec{\mathbf{v}}_m)\Delta\phi_m\|^2&\le \|\nabla(\vec{\mathbf{v}}_m)\|^2\|\Delta\phi_m\|^2_{L^\infty}\\
&\le C\big(\|\nabla\mu_m\cdot\nabla\phi_m\|^2+\|\mu_m\Delta\phi_m\|^2\big)\|\Delta\phi_m\|^2_{H^2}\\
&\le C\big(\|\nabla\mu_m\|_{L^4}^2\|\nabla\phi_m\|^2_{L^4}+\|\mu_m\|_{L^\infty}^2\|\Delta\phi_m\|^2\big)\|\Delta^2\phi_m\|^2\\
&\le C\big(\|\Delta\mu_m\|^2\|\Delta\phi_m\|^2+\|\Delta\mu_m\|^2\|\Delta\phi_m\|^2\big)\|\Delta^2\phi_m\|^2\\
&\le C \|\Delta\mu_m\|^2\|\Delta^2\phi_m\|^2,
\end{aligned}\label{e718}
\end{equation}
where we used the estimate of $\|\Delta\phi_m\|^2$.

The estimate of the third term on the RHS of (\ref{e716}) is more delicate than previous two estimates. First, we have
\begin{equation}
\begin{aligned}
\|\Delta(\vec{\mathbf{v}}_m)\nabla\phi_m\|^2&\le \|\Delta(\mu_m\nabla\phi_m)\|^2\|\nabla\phi_m\|_{L^\infty}^2\\
&\le 4\big(\|\Delta\mu_m\nabla\phi_m\|^2+\|\nabla\mu_m\Delta\phi_m\|^2+\|\mu_m\nabla\Delta\phi_m\|^2\big)\|\nabla\phi_m\|_{L^\infty}^2.
\end{aligned}\label{e719}
\end{equation}

We treat the RHS of (\ref{e719}) term by term. First we have
\begin{equation}
\begin{aligned}
\|\Delta\mu_m\nabla\phi_m\|^2\|\nabla\phi_m\|_{L^\infty}^2&\le \|\Delta\mu_m\|^2\|\nabla\phi_m\|_{L^\infty}^4.
\end{aligned}\label{e720}
\end{equation}
By the Sobolev interpolation inequality $\|f\|_{L^\infty}^2\le C\|f\|_{L^2}\|f\|_{H^3}$ in 3D and (\ref{e321}) we have
\begin{equation}
\begin{aligned}
\|\Delta\mu_m\|^2\|\nabla\phi_m\|_{L^\infty}^4&\le C\|\Delta\mu_m\|^2\|\nabla\phi_m\|_{L^2}^2\|\nabla\phi_m\|_{H^3}^2\\
&\le C\|\Delta\mu_m\|^2\|\nabla\phi_m\|^2\|\Delta^2\phi_m\|^2\\
&\le C\|\Delta\mu_m\|^2\|\Delta^2\phi_m\|^2,
\end{aligned}\label{e721}
\end{equation}
where we have used the estimate of $\|\nabla\phi_m\|^2$.
Therefore, by plugging (\ref{e721}) into (\ref{e720}) we have
\begin{equation}
\|\Delta\mu_m\nabla\phi_m\|^2\|\nabla\phi_m\|_{L^\infty}^2\le C\|\Delta\mu_m\|^2\|\Delta^2\phi_m\|^2.\label{e723}
\end{equation}

For the second and third terms on the RHS of (\ref{e719}), using the same interpolation inequality as above and Lemma 2.2 (X) and (XII) we have
\begin{equation}
\begin{aligned}
\|\nabla\mu_m\Delta\phi_m\|^2\|\nabla\phi_m\|_{L^\infty}^2&\le \|\nabla\mu_m\|_{L^6}^2\|\Delta\phi_m\|_{L^3}^2\|\nabla\phi_m\|_{L^\infty}^2\\
&\le \|\nabla\mu_m\|_{H^1}^2\|\Delta\phi_m\|\|\nabla\Delta\phi_m\|\|\nabla\phi_m\|\|\nabla\phi_m\|_{H^3}\\
&\le \|\Delta\mu_m\|^2\|\nabla\Delta\phi_m\|\|\Delta^2\phi_m\|\\
&\le  \|\Delta\mu_m\|^2\|\Delta^2\phi_m\|^2,
\end{aligned}\label{e724}
\end{equation}
and
\begin{equation}
\begin{aligned}
\|\mu_m\nabla\Delta\phi_m\|^2\|\nabla\phi_m\|_{L^\infty}^2 &\le \|\mu_m\|^2\|\nabla\Delta\phi_m\|_{L^\infty}^2\|\nabla\phi_m\|\|\Delta^2\phi_m\|\\
&\le C\|\nabla\Delta\phi_m\|\|\Delta^3\phi_m\|\|\Delta^2\phi_m\|\\
&\le \frac{1}{96}\|\Delta^3\phi_m\|^2+C\|\nabla\Delta\phi_m\|^2\|\Delta^2\phi_m\|^2.
\end{aligned}\label{e725}
\end{equation}

Plugging (\ref{e723})--(\ref{e725}) into (\ref{e719}) we have
\begin{equation}
\begin{aligned}
\|\Delta(\vec{\mathbf{v}}_m)\nabla\phi_m\|^2\le \frac{1}{24}\|\Delta^3\phi_m\|^2+C\big(\|\nabla\Delta\phi_m\|^2+\|\Delta\mu_m\|^2\big)\|\Delta^2\phi_m\|^2.
\end{aligned}\label{e726}
\end{equation}

Plugging (\ref{e717}), (\ref{e718}) and (\ref{e725}) into (\ref{e716}) we have
\begin{equation}
\|\Delta(\vec{\mathbf{v}}_m\cdot\nabla\phi_m)\|^2\le \frac{1}{6}\|\Delta^3\phi_m\|^2+C\big(\|\nabla\Delta\phi_m\|^2+\|\vec{\mathbf{v}}_m\|^2_{H^1}+\|\Delta\mu_m\|^2\big)\|\Delta^2\phi_m\|^2.\label{e727}
\end{equation}

Plugging (\ref{e727}) into (\ref{e715}), after rearranging terms, we get
\begin{equation}
\begin{aligned}
\frac{d}{dt}\|\Delta^2\phi_m\|^2+\|\Delta^3\phi_m\|^2\le C\big(1+\|\nabla\Delta\phi_m\|^2+\|\vec{\mathbf{v}}_m\|^2_{H^1}+\|\Delta\mu_m\|^2\big)\|\Delta^2\phi_m\|^2.
\end{aligned}\label{e728}
\end{equation}

Combining (\ref{e711}) and (\ref{e728}) we get
\begin{equation}
\begin{aligned}
\frac{d}{dt}\big(\|\nabla\Delta\phi_m\|^2+\|\Delta^2\phi_m\|^2\big)&+\big(\|\nabla\Delta^2\phi_m\|^2+\|\Delta^3\phi_m\|^2\big)\\
&\le C\big(1+\|\nabla\Delta\phi_m\|^2+\|\vec{\mathbf{v}}_m\|^2_{H^1}+\|\Delta\mu_m\|^2\big)\big(\|\nabla\Delta\phi_m\|^2+\|\Delta^2\phi_m\|^2\big).
\end{aligned}\label{e729}
\end{equation}

{\bf Step 5.} Let $\widetilde{T}>0$ be the life span of the solution. For any $0<s\le t\le \widetilde{T}$, applying Gronwall's inequality to (\ref{e729}) over $[s,t]$ and using previous estimates we have
\begin{equation}
\begin{aligned}
&\big(\|\nabla\Delta\phi_m\|^2+\|\Delta^2\phi_m\|^2\big)(t)\\
\le &\exp\bigg\{\int_s^tC\big(1+\|\nabla\Delta\phi_m\|^2+\|\vec{\mathbf{v}}_m\|^2_{H^1}+\|\Delta\mu_m\|^2\big)(\tau)d\tau\bigg\}\big(\|\nabla\Delta\phi_m\|^2
+\|\Delta^2\phi_m\|^2\big)(s)\\
\le &C\big(\|\nabla\Delta\phi_m\|^2
+\|\Delta^2\phi_m\|^2\big)(s).
\end{aligned}\label{e730}
\end{equation}
Integrating (\ref{e730}) w.r.t. $s$ over $(0,t)$ we get
\begin{equation}
\begin{aligned}
t\big(\|\nabla\Delta\phi_m\|^2+\|\Delta^2\phi_m\|^2\big)(t)
\le C\int_0^t\big(\|\nabla\Delta\phi_m\|^2
+\|\Delta^2\phi_m\|^2\big)(s)ds\le C(\widetilde{T}),\ \ \forall\ t\in(0,\widetilde{T}),
\end{aligned}\label{e731}
\end{equation}
which implies that
\begin{equation}
\big(\|\nabla\Delta\phi_m\|^2+\|\Delta^2\phi_m\|^2\big)(t)
\le \frac{C(\widetilde{T})}{t},\ \ \forall\ t\in(0,\widetilde{T}),\label{e732}
\end{equation}
where $C(\widetilde{T})$ is an increasing function of $\widetilde{T}$.

{\bf Step 6.} For any $0<s\le t\le \widetilde{T}$, integrating (\ref{e729}) over time and using (\ref{e732}) we have
\begin{equation}
\begin{aligned}
&\big(\|\nabla\Delta\phi_m\|^2+\|\Delta^2\phi_m\|^2\big)(t)+\int_s^t\big(\|\nabla\Delta^2\phi_m\|^2+
\|\Delta^3\phi_m\|^2\big)(\tau)d\tau\\
\le &\big(\|\nabla\Delta\phi_m\|^2+\|\Delta^2\phi_m\|^2\big)(s)+C(\widetilde{T})\int_s^t\frac{\big(1+
\|\nabla\Delta\phi_m\|^2+\|\vec{\mathbf{v}}_m\|^2_{H^1}+\|\Delta\mu_m\|^2\big)(\tau)}{\tau}d\tau\\
\le &\frac{C(\widetilde{T})}{s}+\frac{C(\widetilde{T})}{s}\int_s^t\big(1+\|\nabla\Delta\phi_m\|^2+
\|\vec{\mathbf{v}}_m\|^2_{H^1}+\|\Delta\mu_m\|^2\big)(\tau)d\tau\le \frac{C(\widetilde{T})}{s}.
\end{aligned}\label{e733}
\end{equation}
In particular, we have
\begin{equation}
\int_s^t\big(\|\nabla\Delta^2\phi_m\|^2+\|\Delta^3\phi_m\|^2\big)(\tau)ds\le \frac{C(\widetilde{T})}{s}.\label{e734}
\end{equation}

For any $\alpha>0$, multiplying (\ref{e734}) by $s^\alpha$ then integrating the result w.r.t. $s$ from $0$ to $t$ we have
\begin{equation}
\int_0^ts^\alpha\int_s^t\big(\|\nabla\Delta^2\phi_m\|^2+\|\Delta^3\phi_m\|^2\big)(\tau)d\tau ds\le C(\widetilde{T})\int_0^ts^{\alpha-1}ds=C(\widetilde{T})\frac{t^\alpha}{\alpha}.\label{e735}
\end{equation}
Changing the order of integration we then have
\begin{equation}
\begin{aligned}
\int_0^ts^\alpha\int_s^t\big(\|\nabla\Delta^2\phi_m\|^2+\|\Delta^3\phi_m\|^2\big)(\tau)d\tau ds&=\int_0^t\bigg(\int_0^\tau s^\alpha ds\bigg)\big(\|\nabla\Delta^2\phi_m\|^2+\|\Delta^3\phi_m\|^2\big)(\tau)d\tau\\
&=\int_0^t\frac{\tau^{\alpha+1}}{\alpha+1}\big(\|\nabla\Delta^2\phi_m\|^2+\|\Delta^3\phi_m\|^2\big)(\tau)d\tau.
\end{aligned}\label{e736}
\end{equation}
Plugging (\ref{e736}) into (\ref{e735}) we have
\begin{equation}
\int_0^t\tau^{\alpha+1}\big(\|\nabla\Delta^2\phi_m\|^2+\|\Delta^3\phi_m\|^2\big)(\tau)d\tau\le C(\widetilde{T})\frac{t^\alpha}{\alpha},\label{e737}
\end{equation}
for any $\alpha>0$ and $0<t<\widetilde{T}$. This completes the proof of Theorem 1.4.

\section{Gevrey regularity in 2D $\&$ 3D}

In this section, we study the Gevrey (spatial) regularity of the solution, and show that the strong solutions obtained in Theorem 1.1--Theorem 1.2 belong to a Gevrey class of regularity defined below. Since the elements in the Gevrey class have high-mode coefficients that decay exponentially in wave number to zero, as a result of Gevrey regularity, one can show that the standard Galerkin method based on the basis (\ref{e21})--(\ref{e23}) converges exponentially fast, see e.g. \cite{DT1,DT2,GST,JMT,LT}. Throughout this section we will assume that $\phi^0(\mathbf{x})\in H^4(\Omega)$. It should be pointed out that in \cite{FT} the authors studied Gevrey regularity for general nonlinear parabolic equations with analytic nonlinearity using the fact that $H^s$ is a Banach algebra for $s>\frac{n}{2}$. But, Foias and Temam \cite{FT-2} get a better result for the Navier-Stokes equations when $u_0\in H^1$, which is not an algebra in either 2D nor 3D, by taking advantage of the quadratic nonlinearity of the NSE. Here, one probably can follow those arguments and reduce the condition on the initial condition. However, this might be very involved since our nonlinearity contains $\phi^3$. So, we follow \cite{FT} instead to simplify the presentation and we do not claim that we have the optimal result regarding the assumption on the initial data. It is also important to note that the eigenfunctions (\ref{e21})--(\ref{e23}) span a subspace of the periodic functions on $[0,2L]\times [0,2l]\times [0,2h]$ which is invariant under the flow generated by (\ref{e11}). Due to this symmetry, the boundary conditions do not explicitly enter the analysis and we can proceed as if (\ref{e11}) were endowed
with periodic boundary conditions.

For the domain $\Omega$ under consideration and $t\ge0$, we define the Gevrey classes $G_t^{p/2}(\Omega)=D(A^{p/2}e^{tA^{1/2}})$, where $A=(I-\Delta)$. These spaces are Hilbert spaces with respect to the inner products
\begin{equation}
(v,w)_{G^{p/2}_t(\Omega)}=\sum_{j\in\mathbf{Z}^n}v_j\cdot \bar{w}_j(1+|j|^2)^pe^{2t(1+|j|^2)^{1/2}}
\end{equation}
with the corresponding norms
\begin{equation}
\|v\|_{G^{p/2}_t(\Omega)}=\bigg(\sum_{j\in\mathbf{Z}^n}|v_j|^2(1+|j|^2)^pe^{2t(1+|j|^2)^{1/2}}\bigg)^{1/2},
\end{equation}
where the $v_j$ and $w_j$ are the Fourier coefficients of $v$ and $w$, respectively, with respect to the basis in (\ref{e21})--(\ref{e23}) where the functions $\cos k\theta$ and $\sin k\theta$ are represented in terms of the exponential function $e^{ik\theta}$ and its conjugate $e^{-ik\theta}$, c.f. \cite{FT}.

We have the following lemma (c.f. \cite{FT}).

\begin{lemma}
For $p\ge2$ and $t\ge0$, $G_t^{p/2}(\Omega)$ is a Banach algebra. That is
$$
\|uv\|_{G^{p/2}_t(\Omega)}\le \|u\|_{G^{p/2}_t(\Omega)}\|v\|_{G^{p/2}_t(\Omega)},
$$
for any $u,v\in {G^{p/2}_t(\Omega)}$.
\end{lemma}

Concerning Gevrey class regularity, we refer the readers to \cite{FT-2} for the solutions of the Navier-Stokes equations and to \cite{FT} for solutions of general nonlinear parabolic equations and to \cite{LT,OT} for the B\'{e}nard convection in a porous medium.

Now we turn to the proof of Theorem 1.5. We first rewrite the equations (\ref{ns}) in terms of $A$ as
\begin{equation}\label{nns}
\left\{
\begin{aligned}
&\frac{d{\phi}_m}{dt}+P_m(\vec{\mathbf{v}}_m\cdot\nabla{\phi}_m)=-A\mu_m+\mu_m,\\
&{\mu}_m=P_m\big[({\phi}_m+\bar{\phi})^3-\overline{({\phi}_m+\bar{\phi})^3}\big]-2{\phi}_m+A{\phi}_m,\\
&\vec{\mathbf{v}}_m=\mathrm{\mathbf{P}}_m\mathrm{\mathbf{P}}_\sigma({\mu}_m\nabla{\phi}_m).
\end{aligned}
\right.
\end{equation}

Applying $A^2e^{tA^{1/2}}$ to $(\ref{nns})_1$ and taking the $L^2$ inner product of the resulting equation with $A^2e^{tA^{1/2}}\phi_m$ we have
\begin{equation}
\begin{aligned}
&\Big(A^2e^{tA^{1/2}}\frac{d\phi_m}{dt},A^2e^{tA^{1/2}}\phi_m\Big)+\big(A^2e^{tA^{1/2}}(\vec{\mathbf{v}}_m\cdot\nabla\phi_m),A^2e^{tA^{1/2}}\phi_m\big)\\
=&-\big(A^2e^{tA^{1/2}}A\mu_m,A^2e^{tA^{1/2}}\phi_m\big)+\big(A^2e^{tA^{1/2}}\mu_m,A^2e^{tA^{1/2}}\phi_m\big).
\end{aligned}\label{e81}
\end{equation}
After integration by parts and using $(\ref{nns})_2$ we have
\begin{equation}
\begin{aligned}
&\frac{1}{2}\frac{d}{dt}\|A^2e^{tA^{1/2}}\phi_m\|^2-\big(A^{5/2}e^{tA^{1/2}}\phi_m,A^2e^{tA^{1/2}}\phi_m\big)
+\big(Ae^{tA^{1/2}}(\vec{\mathbf{v}}_m\cdot\nabla\phi_m),A^{3}e^{tA^{1/2}}\phi_m\big)\\
=&-\big(A^2e^{tA^{1/2}}((\phi_m+\bar{\phi})^3-3\phi_m),A^3e^{tA^{1/2}}\phi_m\big)-\|A^3e^{tA^{1/2}}\phi_m\|^2\\
&+\big(A^2e^{tA^{1/2}}((\phi_m+\bar{\phi})^3),A^2e^{tA^{1/2}}\phi_m\big)-2\|A^2e^{tA^{1/2}}\phi_m\|^2,
\end{aligned}\label{e82}
\end{equation}
which implies that
\begin{equation}
\begin{aligned}
&\frac{1}{2}\frac{d}{dt}\|A^2e^{tA^{1/2}}\phi_m\|^2+\|A^3e^{tA^{1/2}}\phi_m\|^2+2\|A^2e^{tA^{1/2}}\phi_m\|^2\\
\le &\varepsilon \|A^3e^{tA^{1/2}}\phi_m\|^2+C_\varepsilon
\big(\|Ae^{tA^{1/2}}(\vec{\mathbf{v}}_m\cdot\nabla\phi_m)\|^2+\|A^2e^{tA^{1/2}}((\phi_m+\bar{\phi})^3\|^2+\|A^2e^{tA^{1/2}}\phi_m\|^2\big),
\end{aligned}\label{e83}
\end{equation}
where we have used Poincar\'{e}'s inequality to $A^{5/2}e^{tA^{1/2}}\phi_m$. Choosing $\varepsilon=1/2$ we have
\begin{equation}
\begin{aligned}
\frac{d}{dt}\|\phi_m\|^2_{G_t^2}+\|\phi_m\|^2_{G_t^3}+4\|\phi_m\|^2_{G_t^2}
\le C
\big(\|\vec{\mathbf{v}}_m\cdot\nabla\phi_m\|^2_{G_t^1}+\|(\phi_m+\bar{\phi})^3\|^2_{G_t^2}+\|\phi_m\|^2_{G_t^2}\big).
\end{aligned}\label{e84}
\end{equation}

For the second term on the RHS of (\ref{e84}), since $G_t^2$ is a Banach algebra, we have
\begin{equation}
\|(\phi_m+\bar{\phi})^3\|^2_{G_t^2}\le C\|(\phi_m+\bar{\phi})\|^6_{G_t^2}\le C\big(\|\phi_m\|^6_{G_t^2}+|\bar{\phi}|^6 e^{6t}\big).\label{e85}
\end{equation}

Now we deal with the first term on the RHS of (\ref{e84}). Since $G_t^1$ is a Banach algebra, we have
\begin{equation}
\begin{aligned}
\|\vec{\mathbf{v}}_m\cdot A^{1/2}\phi_m\|^2_{G_t^1}&\le \|\vec{\mathbf{v}}_m\|^2_{G_t^1}\|\phi_m\|^2_{G_t^{3/2}}\\
&\le C\|\mu_m A^{1/2}\phi_m\|^2_{G_t^1}\|\phi_m\|^2_{G_t^2}\\
&\le C\|\mu_m\|^2_{G_t^1}\|\phi_m\|^2_{G_t^{3/2}}\|\phi_m\|^2_{G_t^{2}}\\
&\le C\big(\|(\phi_m+\bar{\phi})^3\|^2_{G_t^1}+\|\phi_m\|^2_{G_t^1}+\|\phi_m\|^2_{G_t^2}\big)\|\phi_m\|^4_{G_t^2}\\
&\le C\big(\|\phi_m\|^6_{G_t^1}+|\bar{\phi}|^6 e^{6t}+\|\phi_m\|^2_{G_t^1}+\|\phi_m\|^2_{G_t^2}\big)\|\phi_m\|^4_{G_t^2}\\
&\le C\big(\|\phi_m\|^6_{G_t^2}+|\bar{\phi}|^6 e^{6t}+\|\phi_m\|^2_{G_t^2}\big)\|\phi_m\|^4_{G_t^2},
\end{aligned}\label{e86}
\end{equation}
where we have used Poincar\'{e}'s inequality in various places.

Plugging (\ref{e85}) and (\ref{e86}) into (\ref{e84}) we have
\begin{equation}
\begin{aligned}
\frac{d}{dt}\|\phi_m\|^2_{G_t^2}+\|\phi_m\|^2_{G_t^3}
\le C
\big(\|\phi_m\|^6_{G_t^2}+|\bar{\phi}|^6 e^{6t}+\|\phi_m\|^2_{G_t^2}\big)\big(\|\phi_m\|^4_{G_t^2}+1\big)\equiv P\big(\|\phi_m\|^2_{G_t^2};t\big),
\end{aligned}\label{e87}
\end{equation}
which implies that
\begin{equation}
\begin{aligned}
\|\phi_m(\cdot,t)\|^2_{G_t^2}+\int_0^t\|\phi_m(\cdot,s)\|^2_{G_s^3}ds
\le \int_0^t P\big(\|\phi_m(\cdot,s)\|^2_{G_s^2};s\big)ds+\|\phi_m^0\|^2_{G_0^2}.
\end{aligned}\label{e88}
\end{equation}

Let us consider the first term on the LHS of (\ref{e88}). Since $\|\phi_m^0\|^2_{G_0^2}=\|A^2\phi_m^0\|^2_{L^2}\le C_0\|\phi^0\|^2_{H^4}$, by a continuity argument we know that there exists a time $T_m^{**}>0$ such that
\begin{equation}
\|\phi_m(\cdot,t)\|^2_{G_t^2}\le 2\|\phi_m^0\|^2_{G_0^2}+1\le 2C_0\|\phi^0\|^2_{H^4}+1,\ \  \forall\ t\in[0,T_m^{**}].\label{e89}
\end{equation}
Since $P(\rho;s)$ is an increasing function for $\rho\ge0$ and in $s$, from (\ref{e88}) and (\ref{e89}) we have
\begin{equation}
\|\phi_m(\cdot,t)\|^2_{G_t^2}\le tP\big(2C_0\|\phi^0\|^2_{H^4}+1;t\big)+C_0\|\phi^0\|^2_{H^4},\ \  \forall\ t\in[0,T_m^{**}].\label{e810}
\end{equation}
Let $T^{**}$ be some small value of $t$ such that
\begin{equation}
T^{**}P\big(2C_0\|\phi^0\|^2_{H^4}+1;T^{**}\big)\le {C_0\|\phi^0\|^2_{H^4}+1}.\label{e811}
\end{equation}
Then it holds that
\begin{equation}
\|\phi_m(\cdot,t)\|^2_{G_t^2}\le 2C_0\|\phi^0\|^2_{H^4}+1,\ \  \forall\ t\in[0,T^{**}],\ \ \forall\ m\ge0.\label{e812}
\end{equation}
Upon passing to the limit with $m\to\infty$ we conclude that (\ref{e812}) holds for the limit function. Moreover, from (\ref{e729}) we know that $\phi(t,\mathbf{x})\in H^4(\Omega)$ globally in 2D and locally in 3D if the initial data belong to $H^4(\Omega)$. Therefore, by repeating the process we know that the solution possesses global Gevrey regularity in 2D and local Gevrey regularity in 3D. This completes the proof of Theorem 1.5.

\section{Long-time dynamics}

In this section, we investigate the long-time behavior of the solution to the original problem. First, we show that $\phi$ converges exponentially fast to $\bar{\phi}$ in the $H^1$ norm
as time goes to infinity, provided that $\bar{\phi}$ lies outside the spinodal region where $F''(\phi)=3\phi^2-1\ge0$ (recall $F(\phi)=\frac14(\phi^2-1)^2$) and the initial perturbation $\|\phi^0-\bar{\phi}\|^2_{H^2}$ is sufficiently small.
The proof is in the spirit of \cite{Boyer}. It should be pointed out that the idea for proof has been used in \cite{Boyer} to study the long-time behavior of strong solutions to a coupled Cahn-Hilliard-Navier-Stokes system:
$$
\left\{
\begin{aligned}
&\phi_t+\vec{\mathbf{v}}\cdot\nabla\phi=\Delta \mu,\\
&\mu=F'(\phi)-\Delta\phi,\\
&\vec{\mathbf{v}}_t+(\vec{\mathbf{v}}\cdot\nabla)\vec{\mathbf{v}}+\nabla P=\nu\Delta\vec{\mathbf{v}}+\mu\nabla\phi,\\
&\nabla\cdot \vec{\mathbf{v}}=0,
\end{aligned}
\right.
$$
on bounded domains in 2D and 3D. Similar results are obtained therein, i.e., $\phi$ converges exponentially fast to $\bar{\phi}$ in the $H^1$ norm provided that $\bar{\phi}$ lies outside the spinodal region where $F''(\cdot)\ge0$ and $\|\phi^0-\bar{\phi}\|^2_{H^2}$ is sufficiently small. We briefly explain the idea as follows. Since $\bar{\phi}$ is assumed to be strictly greater than $\sqrt{3}/3$, we note that $F''(\phi)=3\phi^2-1\ge 0$
on a small neighborhood $I_{\bar{\phi}}=[\bar{\phi}-\delta,\bar{\phi}+\delta]$, where $\delta>0$ is sufficiently small. The idea is to solve a modified problem with $F$ replaced by
an auxiliary function $F_{\bar{\phi}}$ whose second order derivative is non-negative on $\mathbb{R}$ and coincides with $F''(\cdot)$ on $I_{\bar{\phi}}$. Then, under the smallness assumption on the initial perturbations, it can be shown that the solution to the modified problem is indeed the solution to the original problem and converges exponentially to the constant state as time goes to infinity. The result holds true for both the 2D and 3D cases.

On the other hand, we show that, in 2D, without the smallness assumption on the initial perturbation, $\phi$ still converges exponentially fast to $\bar{\phi}$ in the $H^2$ norm as time goes to infinity, provided that $L$, the length of the longest edge of the rectangle or the box, is strictly less than $\pi$. This condition will trigger a chain reaction leading the energy estimate performed in Section 4 to a whole new scenario. The result is in strong contrast with the one above. Physically, the result indicates that making the thickness of the diffuse interface relatively large (compared to the dimensions of the domain) leads to constant phase states.

In both situations, the key ingredient of the proofs is the estimate (\ref{e323}) derived in Section 3.

\begin{remark}
In this section, $C$ and $C_i$ will denote generic constants which are independent of the unknown functions and time.
\end{remark}

{\bf Step 1.} We begin with some preparations. First, let us recall
\begin{equation}
F(\phi)=\frac{1}{4}(\phi^2-1)^2.\label{e91}
\end{equation}
Then it is obvious that
\begin{equation}
F'(\phi)=\phi^3-\phi,\ \ F''(\phi)=3\phi^2-1\ \ \mathrm{and\ observe\ that}\ \ \mu=F'(\phi)-\Delta\phi.\label{e92}
\end{equation}

Second, we let
\begin{equation}
\nu=\bar{\phi}=\frac{1}{|\Omega|}\int_\Omega \phi_0(\mathbf{x})d\mathbf{x},\label{e93}
\end{equation}
and assume that
\begin{equation}
|\nu|>\frac{\sqrt{3}}{3}.\label{e94}
\end{equation}

Third, for fixed
\begin{equation}
\delta\in\bigg(0,|\nu|-\frac{\sqrt{3}}{3}\bigg],\label{e95}
\end{equation}
we set
\begin{equation}
I_\delta=[\nu-\delta,\nu+\delta].\label{e96}
\end{equation}
Then it is easy to see that
\begin{equation}
F''(\phi)=3\phi^2-1\ge 0\ \ \mathrm{on}\ I_\delta.\label{e97}
\end{equation}

Fourth, one can easily construct a function $F_\nu$ of $C^3$ class such that
\begin{equation}
F_\nu(\nu)=0,\ \ F'_\nu(\nu)=0,\ \ F''_\nu(z)=F''(z)\ \ \mathrm{on}\ I_\delta,\ \ F''_\nu(z)\ge0\ \ \mathrm{if}\ z\notin I_\delta,\label{e98}
\end{equation}
and such that $F'''_\nu$ is bounded on $\mathbb{R}$. By construction, it is easy to see that
\begin{equation}
F_\nu\ge0.\label{e99}
\end{equation}
By the convexity of $F_\nu$ we see that
\begin{equation}
F_\nu(\nu)-F_\nu(z)\ge F'_\nu(z)(\nu-z),\label{e910}
\end{equation}
which, together with (\ref{e98}), implies that
\begin{equation}
F_\nu(z)\le F'_\nu(z)(z-\nu),\ \ \forall\ z\in\mathbb{R}.\label{e911}
\end{equation}
Moreover, there exist positive constants $F_1$ and $F_2$ such that
\begin{equation}
|F'_\nu(z)|\le F_1+F_2|z|^3,\ \ |F''_\nu(z)|\le F_1+F_2|z|^2,\ \ \forall\ z\in\mathbb{R}.\label{e912}
\end{equation}

{\bf Step 2.} Next, we solve the original problem with $F$ replaced by $F_\nu$ in $(\ref{e11})_2$ and with the same initial data, and denote the solution by $\phi_\nu,\mu_\nu$ and $\vec{\mathbf{v}}_\nu$. In this situation, the estimate (\ref{e311}) is
\begin{equation}
\frac{d}{dt}\bigg(\int_\Omega F_\nu(\phi_\nu)d\mathbf{x}+\frac{1}{2}\|\nabla\phi_\nu\|^2\bigg)+\|\nabla\mu_\nu\|^2+\|\vec{\mathbf{v}}_\nu\|^2
=0.\label{e913}
\end{equation}
Integrating (\ref{e913}) in time we have
\begin{equation}
\begin{aligned}
\int_\Omega F_\nu(\phi_\nu(t))d\mathbf{x}+\frac{1}{2}\|\nabla\phi_\nu(t)\|^2&+\int_0^t\big(\|\nabla\mu_\nu(s)\|^2+\|\vec{\mathbf{v}}_\nu(s)\|^2\big)ds\le \int_\Omega F_\nu(\phi^0)d\mathbf{x}+\frac{1}{2}\|\nabla\phi^0\|^2.
\end{aligned}\label{e914}
\end{equation}
Using (\ref{e911}) and (\ref{e912}) we have
\begin{equation}
\begin{aligned}
\int_\Omega F_\nu(\phi^0)d\mathbf{x}&\le \int_\Omega F'_\nu(\phi^0(\mathbf{x}))(\phi^0(\mathbf{x})-\nu)d\mathbf{x}\\
&\le \|F'_\nu(\phi^0)\|\|\phi^0-\nu\|\\
&\le \big\|F_1+F_2|\phi^0|^3\big\|\|\phi^0-\nu\|\\
&\le C(1+\|\phi^0\|_{L^6}^3)\|\phi^0-\nu\|\\
&\le C(1+\|\phi^0\|_{H^1}^3)\|\phi^0-\nu\|.
\end{aligned}\label{e915}
\end{equation}
Plugging (\ref{e915}) into (\ref{e914}) we then have
\begin{equation}
\begin{aligned}
&\int_\Omega F_\nu(\phi_\nu(t))d\mathbf{x}+\frac{1}{2}\|\nabla\phi_\nu(t)\|^2+\int_0^t\big(\|\nabla\mu_\nu(s)\|^2+\|\vec{\mathbf{v}}_\nu(s)\|^2\big)ds\\
\le &C(1+\|\phi^0\|_{H^1}^3)\|\phi^0-\nu\|+\frac{1}{2}\|\nabla\phi^0\|^2.
\end{aligned}\label{e916}
\end{equation}

{\bf Step 3.} In this step, we exploit the estimate (\ref{e916}) to get more information about the solution. Here we present two different versions of the proof.

{\it Proof 1.} Taking the $L^2$ inner product of $\mu_\nu$ with $\phi_\nu-\nu$ and noting that $(\bar{\mu}_\nu,\phi_\nu-\nu)=0$, where $\bar{\mu}_\nu=\frac{1}{|\Omega|}\int_\Omega \mu_\nu(\mathbf{x},t)d\mathbf{x}$, we have
\begin{equation}
(\mu_\nu-\bar{\mu}_\nu,\phi_\nu-\nu)=\|\nabla\phi_\nu\|^2+\int_\Omega F'_\nu(\phi_\nu)(\phi_\nu-\nu)d\mathbf{x}\ge \|\nabla\phi_\nu\|^2+\int_\Omega F_\nu(\phi_\nu)d\mathbf{x},\label{e917}
\end{equation}
where we have used (\ref{e911}).

Using Poincar\'{e}'s inequality we have
\begin{equation}
(\mu_\nu-\bar{\mu}_\nu,\phi_\nu-\nu)\le \|\mu_\nu-\bar{\mu}_\nu\|\|\phi_\nu-\nu\|\le C\|\nabla\mu_\nu\|\|\nabla\phi_\nu\|\le C\|\nabla\mu_\nu\|^2+ \frac{1}{2}\|\nabla\phi_\nu\|^2,\label{e918}
\end{equation}
which, together with (\ref{e917}), yields
\begin{equation}
\|\nabla\mu_\nu\|^2\ge C\|\nabla\phi_\nu\|^2+C\int_\Omega F_\nu(\phi_\nu)d\mathbf{x}.\label{e919}
\end{equation}

Taking the $L^2$ inner product of $\mu_\nu$ with $-\Delta\phi_\nu$, integrating by parts and using the relevant boundary condition we have
\begin{equation}
(\mu_\nu,-\Delta\phi_\nu)=\|\Delta\phi_\nu\|^2+\int_\Omega F''_\nu(\phi_\nu)|\nabla\phi_\nu|^2d\mathbf{x}\ge \|\Delta\phi_\nu\|^2,\label{e920}
\end{equation}
where we have used (\ref{e98}). By Poincar\'{e}'s inequality we have
\begin{equation}
(\mu_\nu,-\Delta\phi_\nu)=(\nabla\mu_\nu,\nabla\phi_\nu)\le \|\nabla\mu_\nu\|\|\nabla\phi_\nu\|\le C\|\nabla\mu_\nu\|^2+\frac{1}{2}\|\Delta\phi_\nu\|^2.\label{e921}
\end{equation}
Combining (\ref{e920}) with (\ref{e921}) we get
\begin{equation}
\|\nabla\mu_\nu\|^2\ge C\|\Delta\phi_\nu\|^2.\label{e922}
\end{equation}

Therefore, plugging (\ref{e919}) and (\ref{e922}) into (\ref{e916}) we have
\begin{equation}
\begin{aligned}
&\int_\Omega F_\nu(\phi_\nu(t))d\mathbf{x}+\frac{1}{2}\|\nabla\phi_\nu(t)\|^2+\\
&\int_0^t\Big(C\big(\|\nabla\mu_\nu(s)\|^2+\|\Delta\phi_\nu(s)\|^2
+\|\nabla\phi_\nu(s)\|^2+\int_\Omega F_\nu(\phi_\nu(s))d\mathbf{x}\big)+\|\vec{\mathbf{v}}_\nu(s)\|^2\Big)ds\\
\le &C(1+\|\phi^0\|_{H^1}^3)\|\phi^0-\nu\|+\frac{1}{2}\|\nabla\phi^0\|^2,\ \ \forall\ t\ge0.
\end{aligned}\label{e923}
\end{equation}

{\it Proof 2.} Recalling the estimate (\ref{e323}) we have
\begin{equation}
\frac{1}{2}\frac{d}{dt}\|\phi_\nu-\nu\|^2+\|\Delta\phi_\nu\|^2=-\int_\Omega F''(\phi_\nu)|\nabla\phi_\nu|^2d\mathbf{x}.\label{e924}
\end{equation}
Since $F''_\nu\ge0$ on $\mathbb{R}$, we have
\begin{equation}
\frac{1}{2}\frac{d}{dt}\|\phi_\nu-\nu\|^2+\|\Delta\phi_\nu\|^2\le 0,\label{e925}
\end{equation}
which gives
\begin{equation}
\frac{1}{2}\|\phi_\nu(t)-\nu\|^2+\int_0^t\|\Delta\phi_\nu(s)\|^2ds\le \frac{1}{2}\|\phi^0-\nu\|^2.\label{e926}
\end{equation}
Combining (\ref{e916}) and (\ref{e926}) we have
\begin{equation}
\begin{aligned}
&\int_\Omega F_\nu(\phi_\nu(t))d\mathbf{x}+\frac{1}{2}\|\nabla\phi_\nu(t)\|^2+\frac{1}{2}\|\phi_\nu(t)-\nu\|^2+\int_0^t\big(\|\nabla\mu_\nu(s)\|^2
+\|\vec{\mathbf{v}}_\nu(s)\|^2
+\|\Delta\phi_\nu(s)\|^2\big)ds\\
\le &C(1+\|\phi^0\|_{H^1}^3)\|\phi^0-\nu\|+\frac{1}{2}\|\nabla\phi^0\|^2+\frac{1}{2}\|\phi^0-\nu\|^2,\ \ \forall\ t\ge0.
\end{aligned}\label{e927}
\end{equation}

\subsection{Long-time behavior in 2D for small initial perturbations}

We first show the long time dynamics in 2D. Since the solution exists for all time in 2D, when $\|\phi^0-\nu\|_{H^2}^2$ is sufficiently small, thanks to (\ref{e923}) or (\ref{e927}) we have
\begin{equation}
\int_0^t\big(\|\nabla\mu_\nu(s)\|^2+\|\Delta\phi_\nu(s)\|^2\big)ds\le \varepsilon_1,\ \ \forall\ t\ge0.
\label{e928}
\end{equation}
By (\ref{e328}) we have
\begin{equation}
\|\nabla\Delta\phi_\nu\|^2
\le C\big(\|\nabla\mu_\nu\|^2+\|\Delta\phi_\nu\|^2\big),
\label{e929}
\end{equation}
which, together with (\ref{e928}), implies that
\begin{equation}
\int_0^t\|\nabla\Delta\phi_\nu(s)\|^2ds\le \varepsilon_2,\ \ \forall\ t\ge0.\label{e930}
\end{equation}

By (\ref{e412}) we have
\begin{equation}
\frac{1}{2}\frac{d}{dt}\|\Delta\phi_\nu\|^2+\frac{1}{2}\|\Delta^2\phi_\nu\|^2
\le C\|\nabla\Delta\phi_\nu\|^2\|\Delta\phi_\nu\|^2+C\|\nabla\Delta\phi_\nu\|^2,\label{e931}
\end{equation}
where the constant $C$ is independent of time. Applying Gronwall's inequality to (\ref{e931}) we have
\begin{equation}
\|\Delta\phi_\nu(t)\|^2\le \exp\bigg\{C\int_0^t\|\nabla\Delta\phi_\nu(s)\|^2ds\bigg\}\bigg(\|\Delta\phi^0\|^2+C\int_0^t\|\nabla\Delta\phi_\nu\|^2ds\bigg)\le \varepsilon_3,\ \ \forall\ t\ge0,\label{e932}
\end{equation}
where we have used (\ref{e930}) and the smallness of $\|\phi^0-\nu\|_{H^2}^2$. Therefore, by (\ref{e321}) we know that
\begin{equation}
\|\phi_\nu(t)-\nu\|^2_{H^2}\le \varepsilon_4,\ \ \forall\ t\ge0,\label{e933}
\end{equation}
which, together with the Sobolev embedding $H^2\hookrightarrow C^0$, implies that
\begin{equation}
\|\phi_\nu(t)-\nu\|^2_{C^0}\le \varepsilon_5,\ \ \forall\ t\ge0.\label{e934}
\end{equation}

We remark that the constants $\varepsilon_i\ (i=1,\cdots,5)$ are independent of time and go to zero as $\|\phi^0-\nu\|_{H^2}^2$ tends to zero.

We observe that, when $\|\phi^0-\nu\|_{H^2}^2$ is sufficiently small, we can guarantee by (\ref{e95}) that $\phi_\nu\in I_\delta$ for any $t\ge0$ when we choose $\varepsilon_5<\delta$. Thus,
\begin{equation}
F''_\nu(\phi_\nu)=F''(\phi_\nu),\label{e935}
\end{equation}
which, together with (\ref{e98}), implies that
\begin{equation}
F_\nu(\phi_\nu)=F(\phi_\nu)-F(\nu)-(\phi_\nu-\nu) F'(\nu),\ \ \forall\ \phi_\nu\in I_\delta.\label{e936}
\end{equation}
From the definition of the solution we see that adding the above affine function to $F$ does not change the equations. Therefore, we conclude that $\phi_\nu$ is a solution of the original problem. By uniqueness we know that $\phi_\nu=\phi$ for any $t\ge0$ when $\|\phi^0-\nu\|_{H^2}^2$ is sufficiently small.
As a by-product, we know that $\phi\in I_\delta$ for any $t\ge0$.

We note that $\phi$ satisfies the same estimate (\ref{e924}) as $\phi_\nu$ does:
\begin{equation}
\frac{1}{2}\frac{d}{dt}\|\phi-\nu\|^2+\|\Delta\phi\|^2=-\int_\Omega F''(\phi)|\nabla\phi|^2d\mathbf{x}.\label{e937}
\end{equation}
Since $F''\ge0$ on $I_\delta$ and $\phi\in I_\delta$ for any $t\ge0$, we have
\begin{equation}
\frac{1}{2}\frac{d}{dt}\|\phi-\nu\|^2+\|\Delta\phi\|^2\le0,\ \ \forall\ t\ge0.\label{e938}
\end{equation}
Poincar\'{e}'s inequality (c.f. (\ref{e321})) then gives
\begin{equation}
\frac{1}{2}\frac{d}{dt}\|\phi-\nu\|^2+C\|\phi-\nu\|^2\le0,\ \ \forall\ t\ge0,\label{e939}
\end{equation}
from which we derive
\begin{equation}
\|\phi-\nu\|^2(t)\le \|\phi^0-\nu\|^2e^{-Ct},\ \ \forall\ t\ge0.\label{e940}
\end{equation}

Since $\|\phi-\nu\|^2_{H^2}$ is uniformly bounded by virtue of (\ref{e933}), by the interpolation inequality $\|Df\|^2\le C\|D^2f\|\|f\|$ for $\nabla f\cdot\mathbf{n}|_{\partial\Omega}=0$ we have
\begin{equation}
\|\phi-\nu\|^2_{H^1}(t)\le Ce^{-Ct}.\label{e941}
\end{equation}

\subsection{Global existence and long-time behavior in 3D for small initial perturbations}

Now we turn to the long-time behavior of the solution in 3D. First, we need to show that the solution exists for all time when the initial perturbation is sufficiently small. We still work on the auxiliary problem with $F$ replaced by $F_\nu$.

For this purpose, we re-visit the estimate (\ref{e514}), which reads
\begin{equation}
\frac{d}{dt}\|\Delta\phi_\nu\|^2+\|\Delta^2\phi_\nu\|^2
\le C_1\|\Delta\phi_\nu\|^{18}+C_2\|\Delta\phi_\nu\|^2,\label{e942}
\end{equation}
for some constants $C_1,C_2>0$ independent of time. We also have the following estimate from (\ref{e923}) or (\ref{e927}):
\begin{equation}\label{e943}
\int_0^t\|\Delta\phi_\nu(s)\|^2ds\le C_3\|\phi^0-\nu\|_{H^2},\ \ \forall\ t\ge0,
\end{equation}
for some constant $C_3$ which is independent of time. Let
\begin{equation}\label{e944}
\|\phi^0-\nu\|_{H^2}\le \varepsilon,\ \ \mathrm{where}\ \ \varepsilon\le \min\bigg\{\frac{1}{4},\frac{3}{16\big(C_1(\frac{1}{2})^{8}+C_2\big)C_3}\bigg\},
\end{equation}
and
\begin{equation}\label{e945}
T=\sup\big\{\tau\ge0:\ \sup_{0\le t<\tau}\|\Delta\phi_\nu(t)\|^2<1/2\big\}.
\end{equation}
Then it is clear by the continuity argument that $T>0$. If $T=+\infty$, then we have a global solution in 3D. Suppose that $T<+\infty$. Let us focus our discussion over the interval $[0,T)$. From (\ref{e942}) we have
\begin{equation}\label{e946}
\begin{aligned}
\frac{d}{dt}\|\Delta\phi_\nu\|^2+\|\Delta^2\phi_\nu\|^2
&\le \big(C_1\|\Delta\phi_\nu\|^{16}+C_2\big)\|\Delta\phi_\nu\|^2\\
&\le \bigg(C_1\Big(\frac{1}{2}\Big)^{8}+C_2\bigg)\|\Delta\phi_\nu\|^2.
\end{aligned}
\end{equation}
Integrating (\ref{e946}) w.r.t. $t$ we have
\begin{equation}\label{e947}
\begin{aligned}
\|\Delta\phi_\nu(t)\|^2 &\le \bigg(C_1\Big(\frac{1}{2}\Big)^{8}+C_2\bigg)\int_0^t\|\Delta\phi_\nu(s)\|^2ds+\|\Delta\phi^0\|^2\\
&\le \bigg(C_1\Big(\frac{1}{2}\Big)^{8}+C_2\bigg)C_3\|\phi^0-\nu\|_{H^2}+\|\phi^0-\nu\|^2_{H^2}\\
&\le \frac{3}{16}+\frac{1}{16}=\frac{1}{4},\ \ \forall\ t\in[0,T).
\end{aligned}
\end{equation}
This implies that
\begin{equation}\label{e948}
\sup_{0\le t<T}\|\Delta\phi_\nu(t)\|^2\le \frac{1}{4},
\end{equation}
which is a contradiction to the definition of $T$. Therefore, it must holds that $T=+\infty$. We conclude that if the initial data $\phi^0$ satisfies (\ref{e944}), then $\|\Delta\phi_\nu(t)\|^2<1/2$ for all $t\ge0$. In particular, the solution exists globally in time in the 3D case. Moreover, we have from the above
\begin{equation}\label{e949}
\|\Delta\phi_\nu(t)\|^2\le \bigg(C_1\Big(\frac{1}{2}\Big)^{8}+C_2\bigg)C_3\varepsilon+\varepsilon=C_4\varepsilon,\ \ \forall\ t\ge0.
\end{equation}

From (\ref{e949}) and Sobolev embedding we know that when $\varepsilon$ is sufficiently small, $\|\phi_\nu-\nu\|^2_{L^\infty}$ is sufficiently small for all time.
Therefore, $\phi_\nu$ will be confined in the interval $I_\delta$ for all time. Since (\ref{e937}) is valid for 3D also, following the arguments in preceding section we know that
the decay estimate still holds in 3D. In other words, we have
\begin{equation}
\|\phi(t)-\nu\|^2_{H^1}\le Ce^{-Ct},\ \ \forall\ t\ge0.\label{e957}
\end{equation}


\subsection{Long-time behavior in 2D for large initial perturbations}

In this subsection, we study a different aspect of the long-time dynamics of the solution in 2D.
As pointed out in the beginning of this section, we show that, in 2D, without the smallness assumption on the initial perturbation, the solution $\phi$ still converges exponentially fast to $\bar{\phi}$, as time goes to infinity, provided that, as explained in the Introduction, the longest edge of the domain is shorter than $\pi$. Moreover, the decay estimate will be extended to the $H^2$ norm of the solution. We begin with the estimate (\ref{e937}):
\begin{equation}
\begin{aligned}
\frac{1}{2}\frac{d}{dt}\|\phi-\nu\|^2+\|\Delta\phi\|^2&=-\int_\Omega F''(\phi)|\nabla\phi|^2d\mathbf{x}\\
&=-\int_\Omega (3\phi^2-1)|\nabla\phi|^2d\mathbf{x}\\
&\le \|\nabla\phi\|^2\le \frac{L^2}{\pi^2}\|\Delta\phi\|^2,
\end{aligned}\label{e959}
\end{equation}
which gives
\begin{equation}
\frac{1}{2}\frac{d}{dt}\|\phi-\nu\|^2+\Big(1-\frac{L^2}{\pi^2}\Big)\|\Delta\phi\|^2\le0.\label{e960}
\end{equation}
Since $L<\pi$, integrating (\ref{e960}) over time we have
\begin{equation}
\begin{aligned}
\int_0^t\|\Delta\phi(s)\|^2ds\le C,\ \ \forall\ t\ge0.\label{e961}
\end{aligned}
\end{equation}

Since (c.f. (\ref{e313}))
\begin{equation}
\begin{aligned}
\int_0^t\|\nabla\mu(s)\|^2ds\le C,\ \ \forall\ t\ge0,\label{e962}
\end{aligned}
\end{equation}
together with (\ref{e328}), (\ref{e961}) and (\ref{e962}) we have
\begin{equation}
\int_0^t\|\nabla\Delta\phi(s)\|^2ds\le C,\ \ \forall\ t\ge0.\label{e963}
\end{equation}

Recalling (\ref{e412}) we have
\begin{equation}
\frac{1}{2}\frac{d}{dt}\|\Delta\phi\|^2+\frac{1}{2}\|\Delta^2\phi\|^2\le C\|\nabla\Delta\phi\|^2\|\Delta\phi\|^2+C\|\nabla\Delta\phi\|^2.
\label{e964}
\end{equation}
Gronwall's inequality together with (\ref{e961}) and (\ref{e963}) then implies that
\begin{equation}
\|\Delta\phi(t)\|^2+\int_0^t\|\Delta^2\phi(s)\|^2ds\le C,\ \ \forall\ t\ge0.\label{e965}
\end{equation}

Plugging the uniform estimate of $\|\Delta\phi\|^2$ into (\ref{e964}) we have
\begin{equation}
\frac{1}{2}\frac{d}{dt}\|\Delta\phi\|^2+\frac{1}{2}\|\Delta^2\phi\|^2
\le C\|\nabla\Delta\phi\|^2\le C_5\big(\|\nabla\mu\|^2+\|\Delta\phi\|^2\big),\label{e966}
\end{equation}
where we have used (\ref{e328}).

Next, we shall re-visit the estimate (\ref{e311}) which reads
\begin{equation}
\frac{d}{dt}\bigg(\int_\Omega F(\phi)d\mathbf{x}+\frac{1}{2}\|\nabla\phi\|^2\bigg)+\|\nabla\mu\|^2+\|\vec{\mathbf{v}}\|^2
=0.\label{e967}
\end{equation}
We note that since
\begin{equation}
F(\phi)=F(\nu)+F'(\nu)(\phi-\nu)+\frac{1}{2}g(\phi)(\phi-\nu)^2,\ \ \mathrm{where}\ \ g(\phi)=\frac{1}{2}\big(2\nu^2+(\phi+\nu)^2-2\big)\ge -1\label{e968}
\end{equation}
and
\begin{equation}
\int_\Omega (\phi-\nu)d\mathbf{x}=0,\label{e969}
\end{equation}
it holds that
\begin{equation}
\frac{d}{dt}\int_\Omega F(\phi)d\mathbf{x}=\frac{1}{2}\frac{d}{dt}\int_\Omega g(\phi)(\phi-\nu)^2d\mathbf{x}.\label{e970}
\end{equation}
Plugging (\ref{e970}) into (\ref{e967}) we have
\begin{equation}
\frac{d}{dt}\bigg(\frac{1}{2}\int_\Omega g(\phi)(\phi-\nu)^2d\mathbf{x}+\frac{1}{2}\|\nabla\phi\|^2\bigg)+\|\nabla\mu\|^2+\|\vec{\mathbf{v}}\|^2
=0.\label{e971}
\end{equation}

Combining (\ref{e960}) and (\ref{e971}) and dropping the non-negative term $\|\vec{\mathbf{v}}\|^2$ we have
\begin{equation}
\frac{d}{dt}\bigg(\frac{1}{2}\int_\Omega g(\phi)(\phi-\nu)^2d\mathbf{x}+\frac{1}{2}\|\nabla\phi\|^2+\frac{1}{2}\|\phi-\nu\|^2\bigg)+\|\nabla\mu\|^2+(1-L^2/\pi^2)\|\Delta\phi\|^2
\le 0.\label{e972}
\end{equation}
Let
\begin{equation}
E(t)\equiv \frac{1}{2}\int_\Omega g(\phi)(\phi-\nu)^2d\mathbf{x}+\frac{1}{2}\|\nabla\phi\|^2+\frac{1}{2}\|\phi-\nu\|^2\ge -\frac{1}{2}\int_\Omega(\phi-\nu)^2d\mathbf{x}+\frac{1}{2}\|\nabla\phi\|^2+\frac{1}{2}\|\phi-\nu\|^2,\label{e973}
\end{equation}
where we used (\ref{e968}). According to Poincar\'{e}'s inequality we have
\begin{equation}
E(t)\ge \frac{1}{2}\Big(1-\frac{L^2}{\pi^2}\Big)\|\nabla\phi\|^2+\frac{1}{2}\|\phi-\nu\|^2.\label{e974}
\end{equation}

By absorbing the RHS of (\ref{e966}) into the LHS of (\ref{e972}) we have
\begin{equation}
\frac{d}{dt}\bigg(\frac{2C_5}{1-L^2/\pi^2}E(t)+\frac{1}{2}\|\Delta\phi\|^2\bigg)+C_5\big(\|\nabla\mu\|^2+\|\Delta\phi\|^2\big)+\frac{1}{2}\|\Delta^2\phi\|^2
\le 0.\label{e975}
\end{equation}
After integrating in time we have
\begin{equation}
\begin{aligned}
&\bigg(\frac{2C_5}{1-L^2/\pi^2}E(t)+\frac{1}{2}\|\Delta\phi(t)\|^2\bigg)+\int_0^t \bigg(C_5\big(\|\nabla\mu(s)\|^2+\|\Delta\phi(s)\|^2\big)+\frac{1}{2}\|\Delta^2\phi(s)\|^2\bigg) ds\\
\le &\bigg(\frac{2C_5}{1-L^2/\pi^2}E(0)+\frac{1}{2}\|\Delta\phi_0\|^2\bigg).\label{e976}
\end{aligned}
\end{equation}
In particular, we have
\begin{equation}
\|\Delta\phi(t)\|^2\le C,\ \ \forall\ t\ge0,\label{e977}
\end{equation}
which implies that
\begin{equation}
\|\phi(t)-\nu\|^2_{H^2}\le C,\ \ \mathrm{and}\ \ \|\phi(t)-\nu\|_{L^\infty}\le C,\ \ \forall\ t\ge0.\label{e978}
\end{equation}

Thanks to (\ref{e968}) and the definition of the function $g(\phi)$ we have
\begin{equation}
\begin{aligned}
\frac{1}{2}\int_\Omega g(\phi)(\phi-\nu)^2d\mathbf{x}\le &\bigg\|\nu^2+\frac{(\phi+\nu)^2}{2}-1\bigg\|_{L^\infty}\|\phi-\nu\|^2\\
\le &C\|\phi-\nu\|^2,\label{e979}
\end{aligned}
\end{equation}
where we used (\ref{e978}) in the last inequality.
Therefore, by virtue of (\ref{e979}) and (\ref{e974}) we see that
\begin{equation}
E(t)\cong \|\phi-\nu\|^2+\|\nabla\phi\|^2,\label{e980}
\end{equation}
where $\cong$ stands for the equivalence of quantities up to a multiplication by a constant.

Let
\begin{equation}
K(t)\equiv \frac{2C_5}{1-L^2/\pi^2}E(t)+\frac{1}{2}\|\Delta\phi\|^2,\ \ \ Q(t)\equiv C_5\big(\|\nabla\mu\|^2+\|\Delta\phi\|^2\big)+\frac{1}{2}\|\Delta^2\phi\|^2.\label{e981}
\end{equation}
Then from (\ref{e321}) we know that there exists a constant $C_6>0$ such that
\begin{equation}
C_6 K(t)\le Q(t),\label{e982}
\end{equation}
which, together with (\ref{e975}), yields
\begin{equation}
\frac{d}{dt} K(t)+C_6K(t)\le 0.\label{e983}
\end{equation}
Thus, $K(t)$ decays exponentially to zero as time tends to infinity. By (\ref{e981}) we then have
\begin{equation}
\|\phi-\nu\|_{H^2}^2(t)\le Ce^{-Ct}.\label{e984}
\end{equation}

\begin{remark}
The decay estimate obtained in this section can be extended to higher order norms of the solution in 2D, which can not be achieved by the method used in preceding sections. The penalty here is that the size of the domain can not be arbitrary.
\end{remark}

\section{Conclusion}

In this paper, we studied the question of well-posedness and long-time qualitative behavior of a mixture model for solid tumor growth. Specifically, we have shown that: 1) when the initial data belong to $H^2(\Omega)$, strong solutions (see Definition 2.1) to the IBVP are globally (locally resp.) well-posed in 2D (3D resp.); 2) the spatial regularity of the solutions automatically increases by degree 2 within the lifespan of the solutions; 3) when the initial data belong to $H^4(\Omega)$, the solutions are globally (locally resp.) spatially analytic in 2D (3D resp.); 4) in both 2D and 3D, the scalar volume fraction converges exponentially rapidly to its spatial average over the domain, as time goes to infinity, provided that the solution lies outside the spinodal region and the initial perturbation is sufficiently small; and 5) in 2D, if the length of the longest edge of the domain is strictly less than a constant multiple of the interface thickness, then the volume fraction still converges to its spatial average over the domain exponentially rapidly, as time goes to infinity, even if the amplitude of initial perturbation is large. The long-time behavior of the solutions suggest that the distinction between the tumor and the surrounding tissue in the microenvironment will blur and homogenize as time proceeds. On the other hand, it is unknown what happens if the solutions are in the spinodal region. This launches a new interesting problem for future pursue. In addition, we expect to extend the results to the case with non-divergence free velocity fields by adopting similar ideas to the ones presented in this paper.

\

\noindent{\bf Acknowledgement.} The work of E.S. Titi was supported in part by the NSF grants DMS-1009950, DMS-1109640 and DMS-1109645. E.S. Titi also acknowledges the kind hospitality of the Freie
Universit\"at - Berlin, where part of this work was completed; and the support of the Alexander von Humboldt Stiftung/Foundation and the Minerva Stiftung/Foundation. The research of K. Zhao was partially supported by the NSF under agreement No. 0635561. J. Lowengrub gratefully acknowledges partial support from NSF, Div. Math. Sci..

\end{document}